\documentclass[12pt,a4paper,reqno]{article}
\usepackage{geometry}
\usepackage{amsmath}
\usepackage{amsfonts}
\usepackage{theorem}
\usepackage{./documenter_custom}

\usepackage[font=footnotesize,labelfont=bf]{caption}

\usepackage{natbib}

\usepackage{tikz}
\usepackage{ulem}
\usepackage{standalone}
\usetikzlibrary{cd}
\usetikzlibrary{positioning}
\usetikzlibrary{calc}
\usetikzlibrary{fit}
\tikzset{set/.style={draw,circle,inner sep=0pt,align=center}}

\usepackage{xcolor}
\definecolor{morange}{RGB}{255,127,14}
\definecolor{mblue}{RGB}{31,119,180}
\definecolor{mred}{RGB}{214,39,40}
\definecolor{mpurple}{RGB}{148,103,189}
\definecolor{mgreen}{RGB}{44,160,44}

\newtheorem{dfntn}{Definition}

\theoremstyle{remark}
\newtheorem{rmrk}{Remark}

\usepackage[affil-it]{authblk}

\providecommand{\keywords}[1]{\footnotesize\textbf{\textit{Keywords:}} #1}

\makeatletter
\newbox\resumebox
\newenvironment{resume}{
    \global\setbox\resumebox=\vtop\bgroup\fontsize{9}{11}
    \selectfont\advance\hsize-6pc\trivlist\labelsep.5em
    \item[\hskip\labelsep{\fontsize{10}{12}\selectfont\bf R\'esum\'e.}]
    }{
        \endtrivlist\egroup\ifx\@setresume\relax \@setresumea \fi
        }
\def\@setresume{\@setresumea\global\let\@setresume\relax}
\def\@setresumea{\skip@20\p@\advance\skip@-\lastskip\advance\skip@-\baselineskip \vskip\skip@
  \ifvoid\resumebox\else\moveright 3pc \box\resumebox\fi}
\makeatother

\makeatletter
\newbox\abstractbox
\renewenvironment{abstract}{
  \global\setbox\abstractbox=\vtop\bgroup\iflanguage{french}{\selectlanguage{english}}{}
  \fontsize{9}{11}\selectfont 
  \advance \hsize -6pc
  \trivlist 
    \labelsep.5em\item[\hskip\labelsep{\fontsize{10}{12}\selectfont\bf Abstract.}]}
{\endtrivlist\egroup\ifx\@setabstract\relax \@setabstracta \fi}
\def\@setabstract{\@setabstracta\global\let\@setabstract\relax}
\def\@setabstracta{\skip@20\p@ \advance\skip@-\lastskip \advance\skip@-\baselineskip \vskip\skip@
    \moveright 3pc \box\abstractbox}
\makeatother

\makeatletter
\renewcommand\@maketitle{
    {\raggedright 
    \begin{center}
    {\Large \bfseries \sffamily \@title }\\[4ex] 
    {\footnotesize \@author}\\[4ex] 
    \@date\\
    \@setabstract
    \@setresume
    \end{center}} 
}
\makeatother

\theoremstyle{plain}
\newenvironment{acknowledgement}{\par\addvspace{17pt}\small\rmfamily\noindent {\it \ackname}.~}{\par\addvspace{6pt}}
\providecommand{\ackname}{Acknowledgements}

\begin{document}
\normalem
\title{Volume-Preserving Transformers for Learning Time Series Data with Structure}
\author[1]{Benedikt Brantner}\affil[1]{Max-Planck-Institut f\"ur Plasmaphysik, Boltzmannstra\ss{}e 2, 85748 Garching}
\author[2,3]{Guillaume de Romemont}\affil[2]{DAAA, ONERA, Université Paris Saclay, F-92322, Châtillon, France}\affil[3]{Arts et Métiers Institute of Technology, Paris, France}
\author[1]{Michael Kraus}
\author[4]{Zeyuan Li}\affil[4]{Zentrum Mathematik, Technische Universität München, Boltzmannstra\ss{}e 3, 85748 Garching, Germany}
\date{\today}

\maketitle

\begin{abstract} 
    Two of the many trends in neural network research of the past few years have been (i) the learning of dynamical systems, especially with recurrent neural networks such as long short-term memory networks (LSTMs) and (ii) the introduction of transformer neural networks for natural language processing (NLP) tasks. 

While some work has been performed on the intersection of these two trends, those efforts were largely limited to using the vanilla transformer directly without adjusting its architecture for the setting of a physical system.

In this work we develop a transformer-inspired neural network and use it to learn a dynamical system. We (for the first time) change the activation function of the attention layer to imbue the transformer with structure-preserving properties to improve long-term stability. This is shown to be of great advantage when applying the neural network to learning the trajectory of a rigid body.
\end{abstract}
\begin{resume} 
    Deux des nombreuses tendances de la recherche sur les réseaux de neurones de ces dernières années ont été (i) l'apprentissage des systèmes dynamiques, en particulier avec les réseaux de neurones récurrents tels que les `Long Short Term Memory' (LSTM) et (ii) l'introduction de réseaux de neurones de type `transformers' pour le traitement du langage naturel (NLP). 

Bien que certains travaux aient été réalisés à l'intersection de ces deux tendances, ils se sont largement limités à l'utilisation directe du transformer vanilla sans adapter son architecture à la configuration d'un système physique.

Dans ce travail, nous utilisons un réseau de neurones inspiré d'un transformer pour apprendre un système dynamique et, de plus (pour la première fois), nous changeons la fonction d'activation du niveau d'attention afin de lui conférer des propriétés de préservation de la structure dans le but d'améliorer sa stabilité à long terme. Ces propriétés s'avèrent extrêmement importantes lors de l'application du réseau de neurone à la trajectoir d'un corps rigide.
\end{resume}

%

\pagebreak

\section*{Introduction}

\label{16879307210612674453}{}

This work is concerned with the development of accurate and robust neural network architectures for the solution of differential equations. It is motivated by several trends that have determined the focus of research in scientific machine learning [\hyperlinkref{9930491254245683687}{1}] in recent years.

First, machine learning techniques have been successfully applied to the identification of dynamical systems for which data are available, but the underlying differential equation is either (i) not known or (ii) too expensive to solve. The first problem (i) often occurs when dealing with experimental data [\hyperlinkref{3402864322513845913}{2}, \hyperlinkref{8146600536394319917}{3}]; the second one (ii) is crucial in \emph{reduced-order modeling} as will be elaborated on below. Various machine learning models have been shown to be able to capture the behaviour of dynamical systems accurately and (within their regime of validity) make predictions at much lower computational costs than traditional numerical algorithms [\hyperlinkref{9930491254245683687}{1}].

Second, in certain application areas, hitherto established neural network architectures have been gradually replaced by transformer neural networks [\hyperlinkref{2538783340789034789}{4}]. Primarily, this concerns recurrent neural networks such as long short-term memory networks (LSTMs [\hyperlinkref{3150593183825968864}{5}]) that treat time series data, e.g., for natural language processing, but also convolutional neural networks (CNNs) for image recognition [\hyperlinkref{13058814467303700196}{6}]. Transformer networks tend to capture long-range dependencies and contextual information better than other architectures and allow for much more efficient parallelization on modern hardware.

Lastly, the importance of including information about the physical system in scientific machine learning models has been recognized. This did not come as a surprise as it is a finding that had long been established for conventional numerical methods [\hyperlinkref{18037804960457623862}{7}]. In this work, the physical property we consider is \emph{volume preservation} [\hyperlinkref{18037804960457623862}{7}, \hyperlinkref{3006316998437427794}{8}], which is a property of the flow of divergence-free vector fields. There are essentially two approaches through which this property can be accounted for. The first one is the inclusion of terms in the loss function that penalize non-physical behaviour; these neural networks are known as physics-informed neural networks (PINNs\footnotemark{} [\hyperlinkref{3603161085352894112}{9}]). The drawback of this approach is that a given property is enforced only weakly and therefore only approximately satisfied. The other approach, and the one taken here, is to encode physical properties into the network architecture; a related example of this are ``symplectic neural networks" (SympNets [\hyperlinkref{2750591267727325161}{10}]). While this approach is often more challenging to implement, the resulting architectures adhere to a given property exactly.

\footnotetext{See [\hyperlinkref{12908279347149791397}{11}] for an application of PINNs where \emph{non-physical behavior} is penalized.

}

\begin{figure}[h]
\centering
\begin{tikzpicture}[
    squarednodeF/.style={circle, draw=black!60, fill=gray!5, very thick, minimum size=10mm},
    squarednodeR/.style={rectangle, draw=black!60, fill=cyan!5, very thick, minimum size=5mm},
    ]
    \node[squarednodeF]      (full0)                              {$\mathrm{FOM}^0$};
    \node[squarednodeR]        (reduced0)       [below of=full0, yshift=-1cm] {$\mathrm{ROM}^0$};
    \node[squarednodeR]      (reduced1)       [right of=reduced0, xshift=2cm] {$\mathrm{ROM}^1$};
    \node[squarednodeF]        (full1)       [above of=reduced1, yshift=1cm] {$\mathrm{FOM}^1$};
    \node[squarednodeR] (reduced2) [right of=reduced1, xshift=2cm] {$\mathrm{ROM}^2$};
    \node (reduced3) [right of=reduced2, xshift=2cm] {$\cdots$};
    \node (full2) [above of=reduced2, yshift=1cm] {$\cdots$};
    
    \draw[->, mred] (full0.south) to node[right] {} (reduced0.north);
    \draw[->] (reduced0.east) to node[above] {$\mathcal{NN}$} (reduced1.west);
    \draw[->, mred] (reduced1.north) to node[left] {} (full1.south);
    \draw[->, mred] (reduced2.north) to node[left] {} (full2.south);
    \draw[->] (reduced1.east) to node[above] {$\mathcal{NN}$} (reduced2.west);
    \draw[->] (reduced2.east) to node[above] {$\mathcal{NN}$} (reduced3.west);
    \end{tikzpicture}
\caption{Visualization of offline-online split in reduced order modeling. The volume-preserving transformer presented in this work is used for the online phase, i.e. it is used to model $\mathcal{NN}$ in this figure. Here FOM stands for \textit{full-order model} and ROM stands for \textit{reduced-order model}.}
\label{fig:OfflineOnlineSplit}
\end{figure}
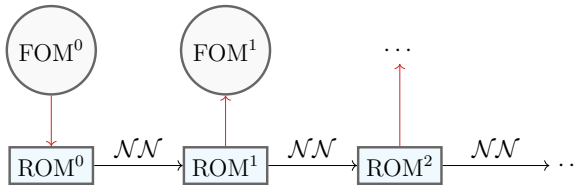

An application, where all of these developments intersect, is \emph{reduced-order modeling} [\hyperlinkref{17010868457598550642}{12}–\hyperlinkref{4173415763474909947}{14}]. Reduced-order modeling is typically split into an \emph{offline phase} and an \emph{online phase} (see \Cref{fig:OfflineOnlineSplit}).

For the online phase we face the following challenges:

\begin{itemize}
\item[1. ] in many cases we need to recover the dynamics of our system from data alone (also known as ``non-intrusive reduced order modeling" [\hyperlinkref{14434722691705629527}{13}–\hyperlinkref{10955802838160861140}{15}]),

\item[2. ] if the big system exhibits specific structure (such as volume-preservation) it is often crucial to also respect this structure in the reduced model [\hyperlinkref{4181908863425119879}{16}, \hyperlinkref{6441307755809283184}{17}].  

\end{itemize}

Our aim in this work is to construct structure-preserving neural network models, that can be used to compute the dynamics of reduced models. Efforts to use neural networks for the online stage have been made in the past, for example using LSTMs [\hyperlinkref{14434722691705629527}{13}, \hyperlinkref{4173415763474909947}{14}].\footnotemark{} Inspired by one of the trends mentioned above, namely the gradual replacement of such architectures, our work will be based on transformers instead.

\footnotetext{Apart from neural networks, there are also other approaches to alleviate the cost during the online stage, notably the ``discrete empirical interpolation method" (DEIM [\hyperlinkref{11021975310350069440}{18}]).

}

While previously, other authors applied transformers for model reduction, and volume-preserving neural networks have been developed as well, both aspects have not yet been considered together. Thus to our knowledge, this work is the first that imbues a transformer with structure-preserving properties\footnotemark{} (namely volume-preservation) and applying it to a system described by a divergence-free vector field.  In the previous work of other authors, transformer neural networks have been used for the online stage of reduced order modeling (i.e. have been applied to dynamical systems) in [\hyperlinkref{11086911267780014317}{19}, \hyperlinkref{4103674955131354433}{20}]. The authors applied the vanilla transformer architecture without taking any physical properties of the system into account. Volume-preserving feedforward neural networks have been developed in [\hyperlinkref{2928453333953959728}{21}]. The authors based the network design on a theorem introduced in [\hyperlinkref{5266324132587918940}{22}] for the design of traditional integrators for divergence-free vector fields [\hyperlinkref{18037804960457623862}{7}].

\footnotetext{In [\hyperlinkref{11823669491487214148}{23}] the transformer is used in a ``structure-conforming" way. In this approach the special architecture of the transformer is leveraged for approximating a solution to a partial differential equation; however no ``structure" is preserved in the way we describe it here, i.e. the neural network-based discretization does not share features with the analytic solution of the differential equation (like volume preservation for example).

}

The remainder of this paper is structured as follows: In \Cref{9901482787949027160} we discuss the basic theory behind divergence-free vector fields and volume-preserving flows, in \Cref{12225886264586945068} the (standard) transformer is introduced, in \Cref{12128329681453367034} we introduce a new class of volume-preserving feedforward neural networks (that differ slightly from what is discussed in e.g. [\hyperlinkref{2928453333953959728}{21}]), in \Cref{1103518901464034041} we introduce our new adapted transformer and in \Cref{13575308147144999289} we finally present results of applying the volume-preserving transformer to a rigid body, an example of a divergence-free vector field (see \Cref{6923796705636635695}). In \Cref{7826571555278476054} we summarize our findings and state potential future investigations to extend the work presented here.

\section{Divergence-free vector fields}

\label{9901482787949027160}{}

The aim of this work is to adapt the transformer architecture to the setting of dynamical systems that are described by \emph{divergence-free vector fields} on \(\mathbb{R}^{d}\), i.e. mappings \(f:\mathbb{R}^d\to\mathbb{R}^d\) for which \(\nabla\cdot{}f = \sum_i\partial_if_i = 0\).  The \emph{flow} of such a vector field \(f\) is the (unique) map \(\varphi_f^t:\mathbb{R}^d\to\mathbb{R}^d\), such that

\begin{equation}
\begin{split}\left. \frac{d}{dt} \right|_{t=t_0}\varphi_f^t(z) = f(\varphi_f^{t_0}(z)),\end{split}\end{equation}

for \(z\in\mathbb{R}^d\) and \(t\) indicates time. For a divergence-free vector field \(f\) the flow \(\varphi_f^t\), in addition to being invertible, is also volume-preserving, i.e. the determinant of the Jacobian matrix of \(\varphi_f^t\) is one: \(\mathrm{det}(D\varphi_f^t) = 1\). This can easily be proved (we drop the subscript \(f\) here):

\begin{equation}
\begin{split}\frac{d}{dt}D\varphi^t(z) = D{}f(\varphi^t(z))D\varphi^t(z) \implies \mathrm{tr}\left( \left(D\varphi^t(z) \right)^{-1} \frac{d}{dt}D\varphi^t(z) \right) = \mathrm{tr}\left(D{}f(\varphi^t(z))\right) = 0.
\end{split}\label{eq:VolumePreservingFlows}\end{equation}

For any matrix-valued function \(A\) we further have

\begin{equation}
\begin{split}\mathrm{tr}(A^{-1}\dot{A}) = \frac{\frac{d}{dt}\mathrm{det}(A)}{\mathrm{det}(A)} \end{split}\end{equation}

and therefore (by using that \(\varphi^t\) is invertible):

\begin{equation}
\begin{split}\frac{d}{dt}\mathrm{det}\left( D\varphi^t(z) \right) = 0.\end{split}\end{equation}

The determinant of \(D\varphi^t\) is therefore constant and we further have \(\det(D\varphi^0) = 1\) because \(\varphi^0\) is the identity. This proves our assertion.

Numerical integrators for ODEs constitute an approximation \(\psi^h\) of the flow \(\varphi_f^t\) where \(h\) denotes the time step, which is fixed in most cases. If the flow \(\varphi_f^t\) exhibits certain properties (like volume preservation) it appears natural to also imbue \(\psi^h\) with these properties. The discipline of doing so is generally known as \emph{geometric numerical integration} [\hyperlinkref{18037804960457623862}{7}, \hyperlinkref{17560367509306694129}{24}].

In recent years, numerical integrators based on neural networks have emerged and it has proven crucial to also imbue these integrators with properties of the system such as symplecticity [\hyperlinkref{2750591267727325161}{10}, \hyperlinkref{7708618804017169771}{25}] and volume preservation [\hyperlinkref{2928453333953959728}{21}]. The neural network architecture presented in \Cref{1103518901464034041} falls in this category.

Let us note that all symplectic and Hamiltonian vector fields\footnotemark{} are also divergence-free but not vice-versa. Symplecticity is a much stronger property than volume preservation. Therefore, preserving symplecticity is often preferable to preserving volume. Still, volume preservation usually offers improved stability and robustness over schemes that do not respect any of the properties of the vector field. Thus, volume-preserving methods can be a viable option when symplectic schemes are not available.

\footnotetext{Strictly speaking Hamiltonian vector fields form a subspace of the space of all symplectic vector fields [\hyperlinkref{12095522061475586953}{26}].

}

\section{The transformer}

\label{12225886264586945068}{}

The transformer architecture [\hyperlinkref{2538783340789034789}{4}] was originally motivated by natural language processing (NLP) tasks and has quickly come to dominate that field. The ``T" in ChatGPT (see e.g. [\hyperlinkref{4243252928265579806}{27}]) stands for ``Transformer" and transformer models are the key element for generative AI. These models are a type of neural network architecture designed to process sequential input data, such as sentences or time-series data. The transformer has replaced, or is in the process of replacing, earlier architectures such as long short-term memory (LSTM) networks [\hyperlinkref{7305568361813525748}{28}] and other recurrent neural networks (RNNs, see [\hyperlinkref{11194799988828072492}{29}]). The transformer architecture is visualized in \Cref{fig:TransformerArchitecture}\footnotemark.

\footnotetext{The three arrows going into the multihead attention module symbolize that the input is used three times: twice when computing the correlation matrix \(C\) and then again when the input is re-weighted based on \(C\). In the NLP literature those inputs are referred to as ``queries", ``keys" and ``values" [\hyperlinkref{2538783340789034789}{4}].

}

As the output of the transformer is of the same dimension as the input to the transformer, we can stack \emph{transformer units} on top of each other; the number of \emph{stacked units} is described by the integer ``L". In essence, the transformer consists of a residual network (ResNet\footnotemark) [\hyperlinkref{5557653725674690011}{30}] and an attention layer. We describe these two core components in some detail.

\footnotetext{ResNets are often used because they improve stability in training [\hyperlinkref{5557653725674690011}{30}] or make it possible to interpret a neural network as an ODE solver [\hyperlinkref{9125255746634034294}{31}].

}

\subsection{Residual Neural Networks}

\label{2534745449722496206}{}

In its simplest form a ResNet is a standard feedforward neural network with an \emph{add connection}:

\begin{equation}
\begin{split}\mathrm{ResNet}: z \rightarrow z + \mathcal{NN}(z),\end{split}\end{equation}

where \(\mathcal{NN}\) is any feedforward neural network. In this work we use a version where the ResNet step is repeated \texttt{n\_blocks} times (also confer \Cref{fig:VolumePreservingFeedForward}), i.e. we have

\begin{equation}
\begin{split}    \mathrm{ResNet} = \mathrm{ResNet}_{\ell_\mathtt{n\_blocks}}\circ\cdots\circ\mathrm{ResNet}_{\ell_2}\circ\mathrm{ResNet}_{\ell_1}.\end{split}\end{equation}

Further, one ResNet layer is simply \(\mathrm{ResNet}_{\ell_i}(z) = z + \mathcal{NN}_{\ell_i}(z) = z + \sigma(W_iz + b_i)\) where we pick tanh as activation function \(\sigma.\)

\subsection{The Attention Layer}

\label{14521051573888138084}{}

The attention layer, which can be seen as a preprocessing step to the ResNet, takes a series of vectors \(z^{(1)}_\mu, \ldots, z^{(T)}_\mu\) as input (the \(\mu\) indicates a specific time sequence) and outputs a \emph{learned convex combination of these vectors}. So for a specific input:

\begin{equation}
\begin{split}input = Z = [z_\mu^{(1)}, z_\mu^{(2)}, \ldots, z_\mu^{(T)}],\end{split}\end{equation}

the output of an attention layer becomes:

\begin{equation}
\begin{split}output = \left[ \sum_{i=1}^Ty^{(1)}_iz_\mu^{(i)}, \sum_{i=1}^Ty^{(2)}_iz_\mu^{(i)}, \ldots, \sum_{i=1}^Ty^{(T)}_iz_\mu^{(i)} \right] ,
\end{split}\label{eq:StandardTransformerOutput}\end{equation}

with the coefficients \(y_{i}^{(j)}\) satisfying \(\sum_{i=1}^Ty^{(j)}_i = 1 \; \forall{} j=1,\ldots,T\). It is important to note that the mapping

\begin{equation}
\begin{split}input \mapsto \left([y_i^{(j)}]_{i=1,\ldots,T,j=1,\ldots,T}\right)\end{split}\end{equation}

is nonlinear. These coefficients are computed based on a correlation of the input data and involve learnable parameters that are adapted to the data during training.

The correlations in the input data are computed through a \emph{correlation matrix} : \(Z \rightarrow Z^T A Z =: C\); they are therefore determined by computing weighted scalar products of all possible combinations of two input vectors where the weighting is done with \(A\); any entry of the matrix \(c_{ij}=(z^{(i)}_\mu)^TAz^{(j)}_\mu\) is the result of computing the scalar product of two input vectors. So any relationship, short-term or long-term, is encoded in this matrix.

After having obtained \(C\), a softmax function is applied column-wise to \(C\) and returns the following output:

\begin{equation}
\begin{split}y_i^{(j)} = [\mathrm{softmax}(C)]_{ij} := e^{c_{ij}}/\left(\sum_{i'=1}^Te^{c_{i'j}}\right).\end{split}\end{equation}

This softmax function maps the correlation matrix to a sequence of \emph{probability vectors}, i.e., vectors in the space \(\mathcal{P}:=\{\mathbf{y}\in[0,1]^d: \sum_{i=1}^dy_i = 1\}\). Every one of these \(d\) probability vectors is then used to compute a convex combination of the input vectors \([z_\mu^{(1)}, z_\mu^{(2)}, \ldots, z_\mu^{(T)}]\), i.e., we get \(\sum_{i=1}^Ty_i^{(j)}z_\mu^{(i)}\) for \(j=1,\ldots,T\). Note that we can also write the convex combination of input vectors as:

\begin{equation}
\begin{split}output = input\Lambda = Z\Lambda,
\end{split}\label{eq:RightMultiplication}\end{equation}

where \(\Lambda = \mathrm{softmax}(C).\) So a linear recombination of input vectors can be seen as a multiplication by a matrix from the right.

Despite its simplicity, the transformer exhibits vast improvements compared to RNNs and LSTMs, including the ability to better capture long-range dependencies and contextual information and its near-perfect parallelizability for computation on GPUs and modern hardware. Furthermore, the simplicity of the transformer architecture makes it possible to interpret all its constituent operations, which is not as easily accomplished when using LSTMs, for example. As the individual operations have a straight-forward mathematical interpretation, it is easier to imbue them with additional structure such as volume-preservation.

\begin{figure}[h]
\centering
\includegraphics[width = .25\textwidth]{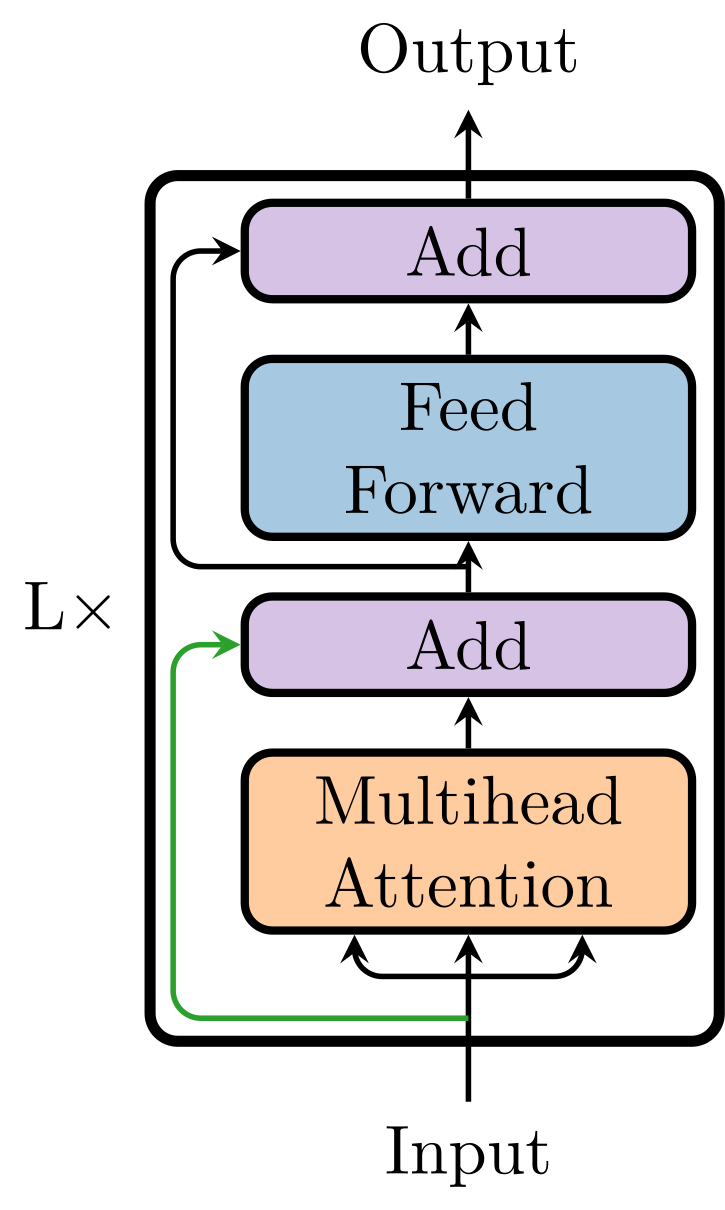}
\caption{Sketch of the transformer architecture. It is a composition of an attention layer and a feedforward neural network. The first \textit{add connection} is drawn in green to emphasize that this can be left out. The Integer ``L'' indicates how often a \textit{transformer unit} (i.e. what is enclosed within the big black borders) is repeated.}
\label{fig:TransformerArchitecture}
\end{figure}

\begin{rmrk} \Cref{fig:TransformerArchitecture} indicates the use of a \emph{multi-head attention layer} as opposed to a \emph{single-head attention layer}. What we described in this section is single-head attention. A multi-head attention layer is slightly more complex: it is a concatenation of multiple single-head attention layers. This is useful for NLP tasks\footnotemark, but introduces additional complexity that makes it harder to imbue the multi-head attention layer with structure-preserving properties. For this reason we stick to single-head attention in this work.\end{rmrk}

\footnotetext{Intuitively, multi-head attention layers allow for attending to different parts of the sequence in different ways (i.e. different heads in the multi-head attention layer \emph{attend to} different parts of the input sequence) and can therefore extract richer contextual information.

}

\section{The Volume-Preserving Transformer}

\label{1103518901464034041}{}

In order to make the transformer volume-preserving we need to modify (i) the ResNet and (ii) the attention layer. We first discuss how we do this for the ResNet and then for the attention layer.

\subsection{Volume-Preserving ResNets}

\label{12128329681453367034}{}

As a first step to constructing a \emph{structure-preserving transformer} we replace the ResNet in the transformer with a feedforward neural network\footnotemark{} that is volume-preserving. The \emph{volume-preserving feedforward layers} here are inspired by the linear and activation modules from [\hyperlinkref{2750591267727325161}{10}]. The key ingredients are upper-triangular matrices \(U\) and lower-triangular matrices \(L\), whose components are such that \(u_{ij} = 0\) if \(i \geq j\) and \(l_{ij} = 0\) if \(i \leq j\), respectively. Similar volume-preserving feedforward neural networks were introduced before [\hyperlinkref{2928453333953959728}{21}]. The difference between the volume-preserving feedforward neural networks in [\hyperlinkref{2928453333953959728}{21}] and the ones presented here is that the ones in [\hyperlinkref{2928453333953959728}{21}] are based on \(G\)-SympNets, whereas ours are based on \(LA\)-SympNets [\hyperlinkref{2750591267727325161}{10}].

\footnotetext{This feedforward neural network is also a ResNet by construction, i.e. the input is again added to the output.

}

Let us consider a \emph{lower-triangular layer}\footnotemark:

\footnotetext{The activation function \(\sigma\) could in principle be chosen arbitrarily, but will be tanh in our experiments.

}

\begin{equation}
\begin{split}x \mapsto x + \sigma(Lx + b) ,
\end{split}\label{eq:VPFF}\end{equation}

with the matrix \(L\) given by:

\begin{equation}
\begin{split}L = \begin{pmatrix}
     0 & 0 & \cdots & 0      \\
     a_{21} & \ddots &        & \vdots \\
     \vdots & \ddots & \ddots & \vdots \\
     a_{n1} & \cdots & a_{n(n-1)}      & 0 
\end{pmatrix}.
\end{split}\label{eq:LinearLower}\end{equation}

The Jacobian matrix of such a layer is of the form

\begin{equation}
\begin{split}J = \begin{pmatrix}
     1 & 0 & \cdots & 0      \\
     b_{21} & \ddots &        & \vdots \\
     \vdots & \ddots & \ddots & \vdots \\
     b_{n1} & \cdots & b_{n(n-1)}      & 1 
\end{pmatrix},\end{split}\end{equation}

and the determinant of \(J\) is \(1\), i.e., the map is volume-preserving.  The same reasoning applies to an ``upper-triangular layer" of the form

\begin{equation}
\begin{split}x \mapsto x + \sigma(Ux + b) .
\end{split}\label{eq:VPFFU}\end{equation}

\begin{figure}
\centering
\includegraphics[width = .31\textwidth]{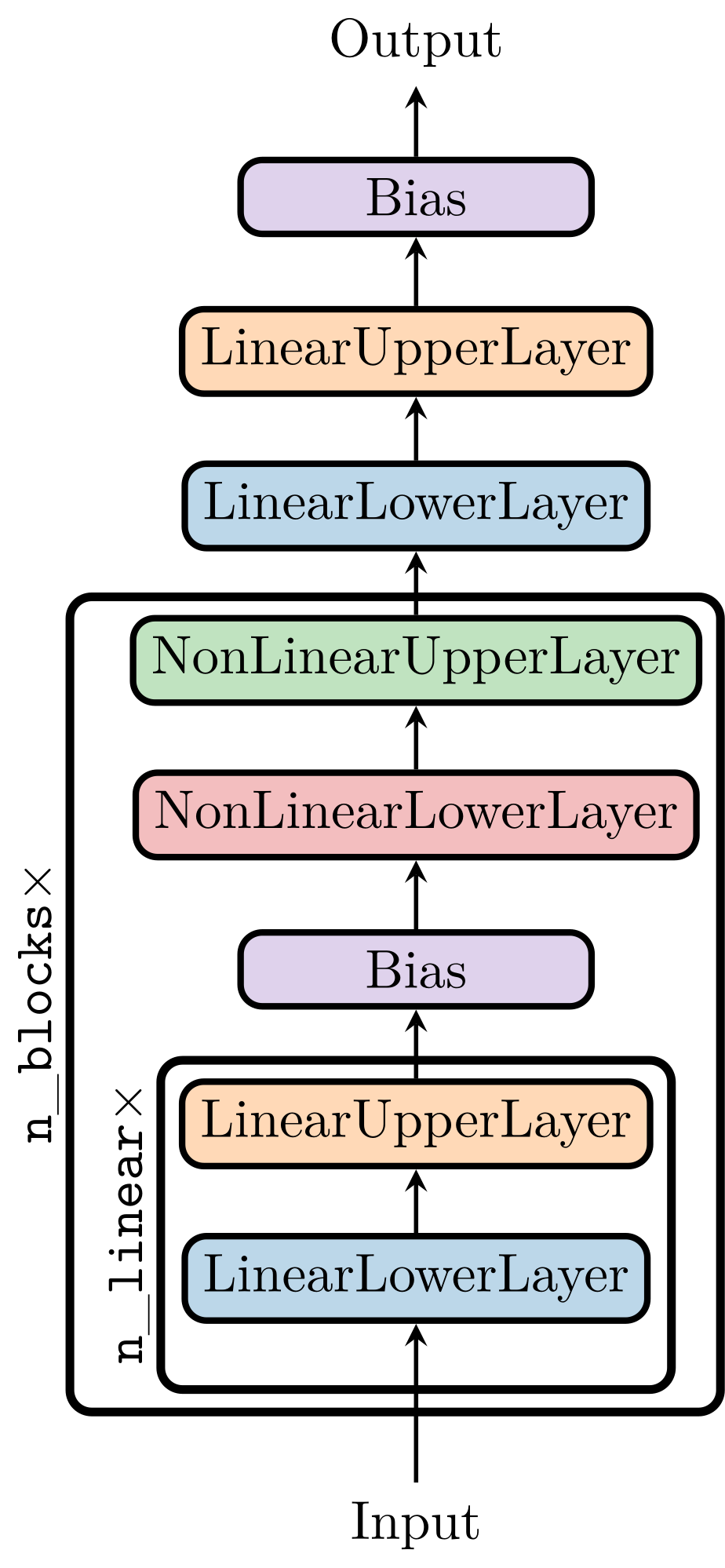}
\caption{Architecture of the volume-preserving feedforward neural network. ``LinearLowerLayer" refers to $x \mapsto x + Lx$ (and similarly for ``LinearUpperLayer"). ``NonLinearLowerLayer" is shown in \Cref{eq:VPFF}. ``Bias" is the addition with a bias vector. Every unit (enclosure within black thick borders) has new parameters for every repetition, so the neural network weights are not repeated.}
\label{fig:VolumePreservingFeedForward}
\end{figure}

In practice we combine many of those layers where the activation function is either (i) a fixed nonlinearity (tanh in our case) or (ii) identity. This is shown in \Cref{fig:VolumePreservingFeedForward}.

Next we introduce the volume-preserving attention layer. We then show how this new layer preserves volume.

\subsection{The Volume-Preserving Attention Layer}

\label{1990617142228263602}{}

The main difficulty in adapting a transformer-like architecture to be volume-preserving is to adapt the activation function. Indeed, the softmax acts vector-wise and cannot preserve volume. We thus replace the softmax by a different activation function, which is based on the Cayley transform:

\begin{equation}
\begin{split}\alpha(Y) = \mathrm{Cayley}(Y) = (\mathbb{I}_{T} - Y)(\mathbb{I}_{T} + Y)^{-1}.\end{split}\end{equation}

The Cayley transform maps skew-symmetric matrices to orthogonal matrices\footnotemark, so if \(Y\) has the property that \(Y^T = -Y,\) then we have \(\alpha(Y)^T\alpha(Y) = \mathbb{I}.\) This can be easily shown:

\footnotetext{The orthogonal matrices \(\{B\in\mathbb{R}^{d\times{}d}:B^TB=\mathbb{I}_d\}\) form a Lie group under regular matrix multiplication. The associated Lie algebra is the vector space of skew-symmetric matrices \(\mathfrak{g}=\{C:C+C^T = \mathbb{O}\}\) and the Lie algebra is mapped to the Lie group via the Cayley transform. More details on this can be found in e.g. [\hyperlinkref{18037804960457623862}{7}].

}

\begin{equation}
\begin{split}\begin{aligned}
\frac{d}{dt}\Big|_{t=0}\alpha(tY)^T\alpha(tY) & = \frac{d}{dt}(\mathbb{I}_{T} - tY)^{-1}(\mathbb{I}_{T} + tY)(\mathbb{I}_{T} - tY)(\mathbb{I}_{T} + tY)^{-1} \\
                                              & = Y + Y - Y - Y = \mathbb{O},
\end{aligned}\end{split}\end{equation}

where we used that \((\Lambda^{-1})^T = (\Lambda^T)^{-1}\) for an arbitrary invertible matrix \(\Lambda\). We now define a new activation function for the attention mechanism which we denote by

\begin{equation}
\begin{split}\Lambda(Z) = \alpha (Z^T A Z),
\end{split}\label{eq:VolumePreservingActivation}\end{equation}

where \(A\) is a \emph{learnable skew-symmetric matrix}. Note that the input into the Cayley transform has to be a skew-symmetric matrix in order for the output to be orthogonal; hence we have the restriction on \(A\) to be skew-symmetric. With this, \(\Lambda(Z)\) is an orthogonal matrix and the entire mapping is equivalent to a multiplication by an orthogonal matrix in the \emph{vector representation} shown in \Cref{eq:isomorphism}. Note that the attention layer can again be seen as a reweighting of the input sequence and thus as a multiplication by a matrix from the right (see \Cref{eq:RightMultiplication} and the comment below that equation):

\begin{equation}
\begin{split}Z \mapsto Z\Lambda(Z).
\end{split}\label{eq:LambdaRight}\end{equation}

We conclude by formalizing what the volume-preserving attention mechanism:

\begin{dfntn} The \textbf{volume-preserving attention mechanism} is a map based on the Cayley transform that reweights a collection of input vectors \(Z\) via \(Z \mapsto Z\Lambda(Z),\) where \(\Lambda(Z) = \mathrm{Cayley}(Z^TAZ)\) and \(A\) is a learnable skew-symmetric matrix.  \end{dfntn}

\subsection{How is Structure Preserved?}

\label{18309819361304734052}{}

Here we discuss how \emph{volume-preserving attention} preserves structure. To this end, we first have to define what volume preservation means for the product space

\begin{equation}
\begin{split}\times_T \mathbb{R}^{d} \equiv \underbrace{ \mathbb{R}^{d} \times \cdots \times \mathbb{R}^{d} }_{\text{$T$ times}} .\end{split}\end{equation}

Note that the input to the transformer \(Z\) is an element of this product space. Now consider the isomorphism \(\hat{}: \times_T \mathbb{R}^{d} \stackrel{\approx}{\longrightarrow} \mathbb{R}^{dT}\) of the form:

\begin{equation}
\begin{split}Z = \begin{pmatrix}
            z_1^{(1)} &  z_1^{(2)} & \quad\cdots\quad & z_1^{(T)} \\
            z_2^{(1)} &  z_2^{(2)} & \cdots & z_2^{(T)} \\
            \cdots &  \cdots & \cdots & \cdots \\
            z_d^{(1)} & z_d^{(2)} & \cdots & z_d^{(T)}
    \end{pmatrix}
\mapsto \hat{Z} = 
\begin{bmatrix}
    z_1^{(1)} \\
    z_1^{(2)} \\
    \cdots \\
    z_1^{(T)} \\
    z_2^{(1)} \\
    \cdots \\
    z_d^{(T)}
\end{bmatrix} 
=: Z_\mathrm{vec}.
\end{split}\label{eq:isomorphism}\end{equation}

We refer to the inverse of \(Z \mapsto \hat{Z}\) as \(Y \mapsto \check{Y}\). In the following we also write \(\hat{\varphi}\) for the mapping \(\,\hat{}\circ\varphi\circ\check{}\,\).

\begin{dfntn} A mapping \(\varphi: \times_T \mathbb{R}^{d} \to \times_T \mathbb{R}^{d}\) is said to be \textbf{volume-preserving} if the associated \(\hat{\varphi}\) is volume-preserving. \end{dfntn}

In the transformed coordinate system, that is in terms of the vector \(Z_\mathrm{vec}\) defined in \Cref{eq:isomorphism}, this is equivalent to multiplication by a block-diagonal matrix \(\widehat{\Lambda(Z)}\) from the left:

\begin{equation}
\begin{split}\widehat{\Lambda(Z)} Z_\mathrm{vec} :=
\begin{pmatrix}
\Lambda(Z)^T & \mathbb{O} & \cdots  & \mathbb{O} \\
\mathbb{O} & \Lambda(Z)^T & \cdots & \mathbb{O} \\
\cdots & \cdots & \ddots & \cdots \\ 
\mathbb{O} & \mathbb{O} & \cdots & \Lambda(Z)^T \\
\end{pmatrix}
\left[\begin{array}{c}  z_1^{(1)} \\ z_1^{(2)} \\ \ldots \\ z_1^{(T)} \\ z_2^{(1)} \\ \ldots \\ z_d^{(T)} \end{array}\right] .
\end{split}\label{eq:LambdaApplication}\end{equation}

and it is easy to see that \(\widehat{\Lambda(Z)}\) in \Cref{eq:LambdaApplication} is an orthogonal matrix. We show that \Cref{eq:LambdaApplication} is true:

\begin{equation}
\begin{split}    \widehat{Z\Lambda(Z)} \equiv \widehat{\sum_{k=1}^Tz_i^{(k)}\lambda_{kj}} \equiv \begin{bmatrix} \sum_{k=1}^T z_1^{(k)}\lambda_{k1} \\ \sum_{k=1}^T z_1^{(k)}\lambda_{k2} \\ \ldots \\ \sum_{k=1}^T z_1^{(k)}\lambda_{kT} \\ \sum_{k=1}^T z_2^{(k)}\lambda_{k1} \\ \ldots \\ \sum_{k=1}^T z_d^{(k)}\lambda_{kT} \end{bmatrix} = \begin{bmatrix} \sum_{k=1}^T \lambda_{k1}z_1^{(k)} \\ \sum_{k=1}^T \lambda_{k2}z_1^{(k)} \\ \ldots \\ \sum_{k=1}^T \lambda_{kT}z_1^{(k)} \\ \sum_{k=1}^T \lambda_{k1}z_2^{(k)} \\ \ldots \\ \sum_{k=1}^T \lambda_{kT}z_d^{(k)} \end{bmatrix} = \begin{bmatrix} [\Lambda(Z)^T z_1^{(\bullet)}]_1 \\ [\Lambda(Z)^T z_1^{(\bullet)}]_2 \\ \ldots \\ [\Lambda(Z)^T z_1^{(\bullet)}]_T \\ [\Lambda(Z)^T z_2^{(\bullet)}]_1 \\ \ldots \\ [\Lambda(Z)^T z_d^{(\bullet)}]_T \end{bmatrix} = \begin{bmatrix} \Lambda(Z)^Tz_1^{(\bullet)} \\ \Lambda(Z)^Tz_2^{(\bullet)} \\ \ldots \\ \Lambda(Z)^Tz_d^{(\bullet)} \end{bmatrix},
    \end{split}\label{eq:ProductStructure}\end{equation}

where we defined:

\begin{equation}
\begin{split}    z_i^{(\bullet)} := \begin{bmatrix} z_i^{(1)} \\ z_i^{(2)} \\ \ldots \\ z_i^{(T)} \end{bmatrix}.\end{split}\end{equation}

\begin{figure}
\centering
\includegraphics[width = .29\textwidth]{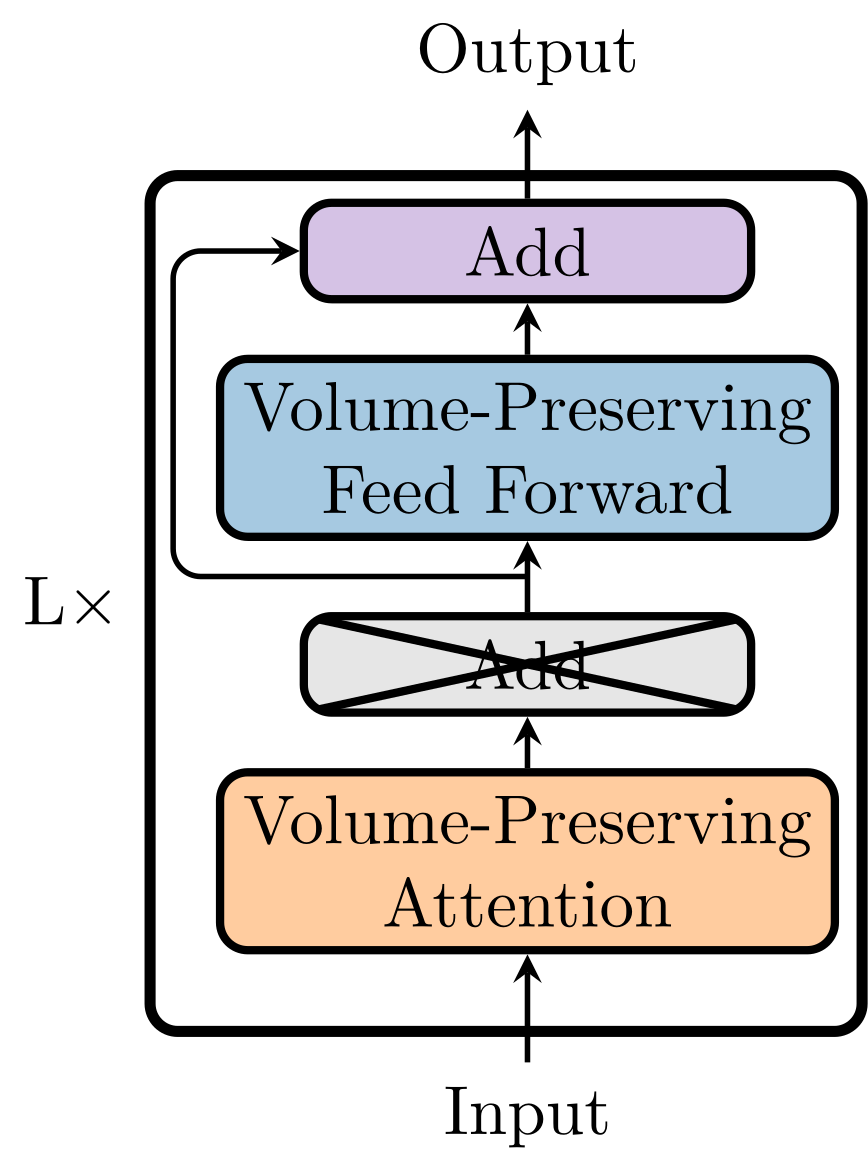}
\caption{Architecture of the volume-preserving transformer. In comparison with the standard transformer in~\Cref{fig:TransformerArchitecture}, (i) the feedforward layer has been replaced with the volume-preserving feedforward neural network from~\Cref{fig:VolumePreservingFeedForward}, (ii) the attention layer has been replaced with a volume-preserving attention layer and (iii) the Add layer has been removed. Similar to~\Cref{fig:TransformerArchitecture} the integer ``L'' indicates how often a \textit{transformer unit} is repeated.}   
\label{fig:VolumePreservingTransformerArchitecture}
\end{figure}

Also note that in \Cref{eq:ProductStructure} the expression after the first ``\(\equiv\)" sign the \((i,j)\)-th element of the matrix, not the entire matrix (unlike the other terms in \Cref{eq:ProductStructure}).

\begin{rmrk} To our knowledge there is no literature on volume-preserving multi-step methods. There is however significant work on \emph{symplectic multi-step methods} [\hyperlinkref{18037804960457623862}{7}, \hyperlinkref{10009329608770522327}{32}, \hyperlinkref{9298822944544965621}{33}]. Of the two definitions of symplecticity for multi-step methods given in [\hyperlinkref{18037804960457623862}{7}], that of \emph{\(G\)-symplecticity} is similar to the definition of volume preservation given here as it is also defined on a product space. The product structure through which we defined volume preservation also bears strong similarities to ``discrete multi-symplectic structures" defined in [\hyperlinkref{10955802838160861140}{15}, \hyperlinkref{9074255080583970709}{34}]. \end{rmrk}

The main result of this section was to show that this attention mechanism is \emph{volume-preserving} in a product space that is spanned by the input vectors. We made the transformer volume-preserving  by (i) replacing the feedforward neural network (ResNet) by a volume-preserving feedforward network (volume-preserving ResNet), (ii) replacing standard attention by volume-preserving attention and (iii) removing the first add connection\footnotemark. The resulting architecture is shown in \Cref{fig:VolumePreservingTransformerArchitecture}.

\footnotetext{Removal of the add connection is necessary as the addition with the input is not a volume-preserving operation. 

}

\section{Experimental results}

\label{13575308147144999289}{}

\begin{figure}[h]
\centering
\includegraphics[width=.65\textwidth]{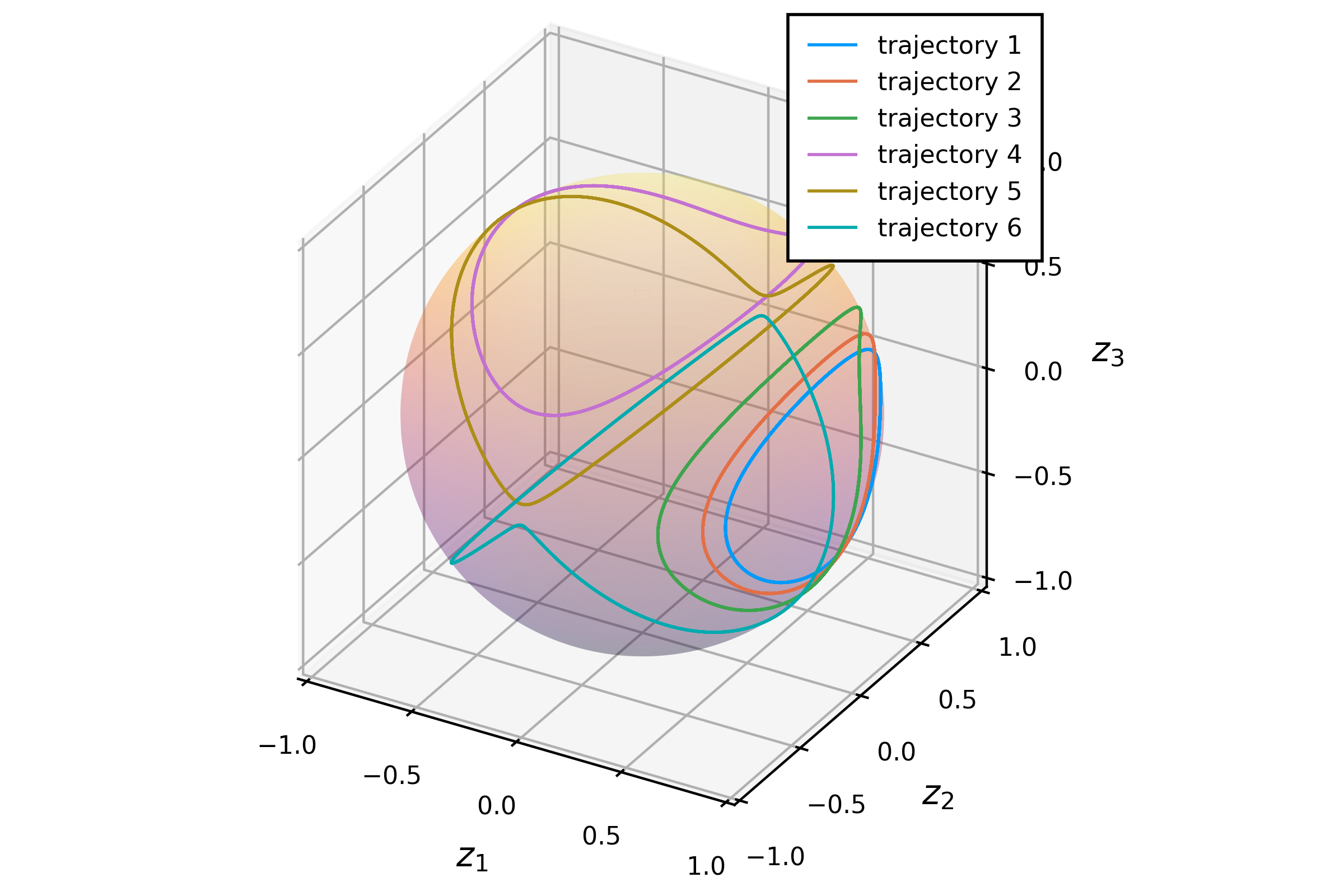}
\caption{Rigid body trajectories for $\mathfrak{a} = 1$, $\mathfrak{b} = -1/2$ and $\mathfrak{c} = -1/2$ and different initial conditions.}
\label{fig:RigidBodyCurves}
\end{figure}

In the following, we will consider the rigid body as an example to study the performance of our new volume-preserving transformer. We will solve the following equations (see \Cref{6923796705636635695} for the derivation):

\begin{equation}
\begin{split}\frac{d}{dt}\begin{bmatrix} z_1 \\  z_2 \\ z_3  \end{bmatrix} 
= \begin{bmatrix} \mathfrak{a}z_2z_3 \\ \mathfrak{b}z_1z_3 \\ \mathfrak{c}z_1z_2 \end{bmatrix} ,
\end{split}\label{eq:RigidBodyEquations}\end{equation}

with \(\mathfrak{a} = 1\), \(\mathfrak{b} = -1/2\) and \(\mathfrak{c} = -1/2\). We immediately see that the vector field in \Cref{eq:RigidBodyEquations} is trivially divergence-free. In \Cref{fig:RigidBodyCurves} we show some trajectories.

We will compare three different neural network architectures that are trained on simulation data of the rigid body.  These architectures are:

\begin{table}[h]
\centering
\footnotesize\begin{tabulary}{\linewidth}{L C C C C}
\toprule
Architecture & \texttt{n\_linear} & \texttt{n\_blocks} & L & Total number of parameters \\
\toprule
Volume-preserving feedforward & 1 & 6 & - & 135 \\
Volume-preserving transformer & 1 & 2 & 3 & 162 \\
Standard transformer & - & 2 & 3 & 213 \\
\bottomrule
\end{tabulary}

\end{table}

\begin{rmrk} Using transformer neural networks instead of standard feedforward neural networks for \emph{integrating} ordinary differential equations can be motivated similarly to using multi-step methods as opposed to single-step methods in traditional numerical integration (another motivation comes from the possibility to consider parameter-dependent equations as discussed below). In [\hyperlinkref{15842744304998993410}{35}] it is stated that multi-step methods constitute ``[an entire class] of integration algorithms of arbitrary order of approximation accuracy that require only a single function evaluation in every new step". We conjecture that this also holds true for transformer-based integrators: we can hope to build higher-order methods without increased cost. \end{rmrk}

For the standard transformer, we further remove the optional add connection (i.e. the green line in \Cref{fig:TransformerArchitecture}) to have a better comparison with the volume-preserving transformer which does not have an add connection. For the standard transformer, \texttt{n\_blocks} refers to the number of ResNet layers we use (the last ResNet layer always has a linear activation). The activation functions in the \emph{feedforward layer} (see \Cref{fig:TransformerArchitecture}) and volume-preserving feedfoward layers (the \emph{non-linear layers} in \Cref{fig:VolumePreservingFeedForward} and the \emph{volume-preserving feedforward layers} in \Cref{fig:VolumePreservingTransformerArchitecture}) are all tanh. For the standard transformer and the volume-preserving transformer we further pick \(T = 3\), i.e. we always feed three time steps into the network during training and validation. We also note that strictly speaking \(T\) is not a hyperparameter of the network as its choice does not change the architecture: the dimensions of the matrix \(A\) in the volume-preserving activation in \Cref{eq:VolumePreservingActivation} (or the equivalent for the standard attention mechanism) are independent of the number of time steps \(T\) that we feed into the network.

For nondegenerate matrices that are of size \((1\times1)\) to \((5\times5)\), or in our notation for \(T = 1, \ldots, 5\), we use explicit matrix inverses. We here state that such an explicit inverse always exists by using matrix adjugates. A proof can be found in e.g. [\hyperlinkref{16425831850931499293}{36}, Proposition VIII.4.16]. For a (\(1\times1\)) matrix this inverse is simply:

\begin{equation}
\begin{split}a \mapsto a^{-1},\end{split}\end{equation}

and for a (\(2\times2\)) matrix it is:

\begin{equation}
\begin{split}\begin{pmatrix} a & b \\ c & d \end{pmatrix}^{-1} = \frac{1}{ad - bc}\begin{pmatrix} d & -b \\ -c & a \end{pmatrix}.
\end{split}\label{eq:inverse2}\end{equation}

For matrices of increasing size this explicit expression gets increasingly expensive, which is why in \texttt{GeometricMachineLearning.jl} we perform it explicitly only for matrices up to size \(5\times5\) (we generate these explicit expressions with \texttt{Symbolics.jl} [\hyperlinkref{7730444816824140145}{37}]). For inverting bigger matrices we use LU decompositions [\hyperlinkref{10297842816248433659}{38}]. We also note that such an explicit inverse is easily parallelizable on GPU as \Cref{eq:inverse2} (and other inverses) already constitutes the computations performed in a GPU kernel.

\subsection{Training data}

\label{15136827737132712471}{}

As training data we take solutions of \Cref{eq:RigidBodyEquations} for various initial conditions:

\begin{equation}
\begin{split}\mathtt{ics} = \left\{ \begin{bmatrix} \sin(v) \\ 0 \\ \cos(v) \end{bmatrix}, \begin{bmatrix} 0 \\ \sin(v) \\ \cos(v) \end{bmatrix}: v\in0.1:0.01:2\pi \right\},
\end{split}\label{eq:Ics}\end{equation}

where \(v\in0.1:0.01:2\pi\) means that we incrementally increase \(v\) from 0.1 to \(2\pi\) by steps of size 0.01. We then integrate \Cref{eq:RigidBodyEquations} for the various initial conditions in \Cref{eq:Ics} with implicit midpoint for the interval \([0,12]\) and a step size of \(0.2\). The integration is done with \texttt{GeometricIntegrators.jl} [\hyperlinkref{14370501854957973096}{39}]. In \Cref{fig:RigidBodyCurves}, we show some of the curves for the following initial conditions:

\begin{equation}
\begin{split}\left\{
\begin{bmatrix} \sin(1.1) \\  0.       \\  \cos(1.1)\end{bmatrix},
\begin{bmatrix} \sin(2.1) \\  0.       \\  \cos(2.1)\end{bmatrix},
\begin{bmatrix} \sin(2.2) \\  0.       \\  \cos(2.2)\end{bmatrix},
\begin{bmatrix}  0.       \\ \sin(1.1) \\  \cos(1.1)\end{bmatrix},
\begin{bmatrix}  0.       \\ \sin(1.5) \\  \cos(1.5)\end{bmatrix}, 
\begin{bmatrix}  0.       \\ \sin(1.6) \\  \cos(1.6)\end{bmatrix}
\right\}.\end{split}\end{equation}

All solutions lie on a sphere of radius one. That they do is a special property of the rigid body (see \Cref{6923796705636635695}) equation and is proofed in [\hyperlinkref{18037804960457623862}{7}, Theorem IV.1.6] for example.

\subsection{Loss functions}

\label{1363865436639116350}{}

For training the feedforward neural network and the transformer we pick similar loss functions. In both cases they take the form:

\begin{equation}
\begin{split}L_{\mathcal{NN}}(input, output) = \frac{ ||output - \mathcal{NN}(input)||_2 }{ ||output||_2 },\end{split}\end{equation}

where \(||\cdot||_2\) is the \(L_2\)-norm. The only difference between the two losses for the feedforward neural network and the transformer is that \(input\) and \(output\) are vectors \(\in\mathbb{R}^d\) in the first case and matrices \(\in\mathbb{R}^{d\times{}T}\) in the second.

\subsection{Details on Training and Choice of Hyperparameters}

\label{5345627476819079035}{}

The code is implemented in Julia [\hyperlinkref{287835226255993214}{40}] as part of the library \texttt{GeometricMachineLearning.jl} [\hyperlinkref{17860356218777098896}{41}]. All the computations performed here are done in single precision on an NVIDIA Geforce RTX 4090 GPU [\hyperlinkref{13147958005838726490}{42}]. We use \texttt{CUDA.jl} [\hyperlinkref{14686158125302097530}{43}] to perform computations on the GPU.

We train the three networks for \(5\cdot10^5\) epochs and use an Adam optimizer [\hyperlinkref{6288960670092762841}{44}] with adaptive learning rate \(\eta\):

\begin{equation}
\begin{split}\eta = \exp\left(log\left(\frac{\eta_2}{\eta_1}\right) / \mathtt{n\_epochs}\right)^t\eta_1,\end{split}\end{equation}

where \(\eta_1\) is the initial learning rate and \(\eta_2\) is the final learning rate. We used the following values for the hyperparameters (mostly taken from [\hyperlinkref{12925252812113732274}{45}]):

\begin{table}[h]
\centering
\footnotesize\begin{tabulary}{\linewidth}{R L L L L L L}
\toprule
Name & \(\eta_1\) & \(\eta_2\) & \(\rho_1\) & \(\rho_2\) & \(\delta\) & \texttt{n\_epochs} \\
\toprule
Value & \(10^{-2}\) & \(10^{-6}\) & \(0.9\) & \(0.99\) & \(10^{-8}\) & \(5\cdot10^5\) \\
\bottomrule
\end{tabulary}

\end{table}

With these settings we get the following training times (given as HOURS:MINUTES:SECONDS) for the different networks:

\begin{table}[h]
\centering
\footnotesize\begin{tabulary}{\linewidth}{R L L L}
\toprule
Architecture & VPFF & VPT & ST \\
\toprule
Training time & 4:02:09 & 5:58:57 & 3:58:06 \\
\bottomrule
\end{tabulary}

\end{table}

\begin{figure}
\centering
\includegraphics[width = .6\textwidth]{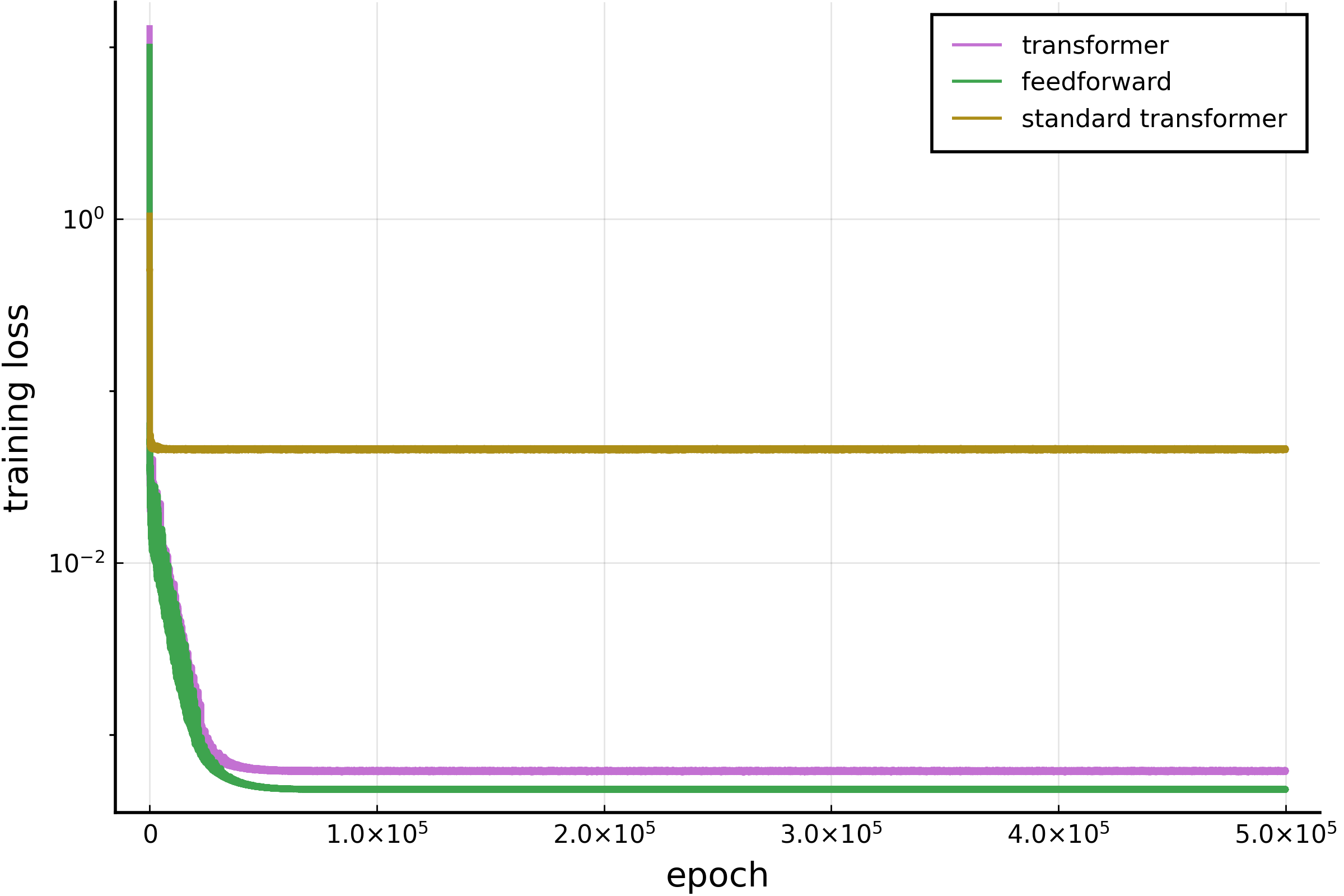}
\caption{Training loss for the different networks.}
\label{fig:TrainingLoss}
\end{figure}

The time evolution of the different training losses is shown in \Cref{fig:TrainingLoss}. The training losses for the volume-preserving transformer and the volume-preserving feedforward neural network reach very low levels (about \(5 \times 10^{-4}\)), whereas the standard transformer is stuck at a rather high level (\(5 \times 10^{-2}\)). In addition to the Adam optimizer we also tried stochastic gradient descent (with and without momentum) and the BFGS optimizer [\hyperlinkref{16796601878836832165}{46}], but obtained the best results with the Adam optimizer. We also see that training the VPT takes longer than the ST even though it has fewer parameters. This is probably because the softmax activation function requires fewer floating-point operations than the inverse in the Cayley transform. We hence observe that even though \texttt{GeometricMachineLearning.jl} has an efficient explicit matrix inverse implemented, this is still slower than computing the exponential in the softmax.

\subsection{Using the Networks for Prediction}

\label{8603635818716639636}{}

\begin{figure}
\includegraphics[width = .33\textwidth]{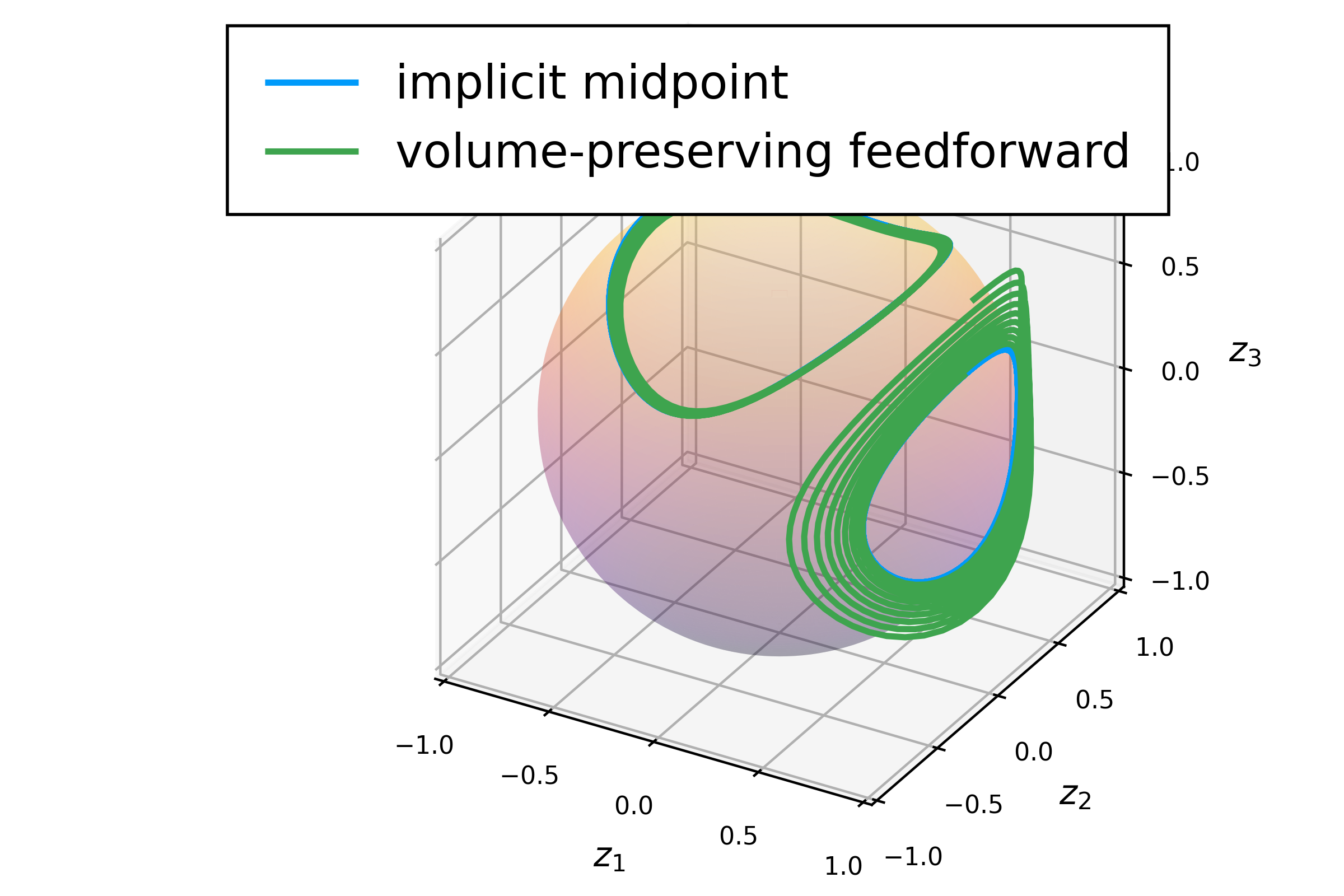}%
\includegraphics[width = .33\textwidth]{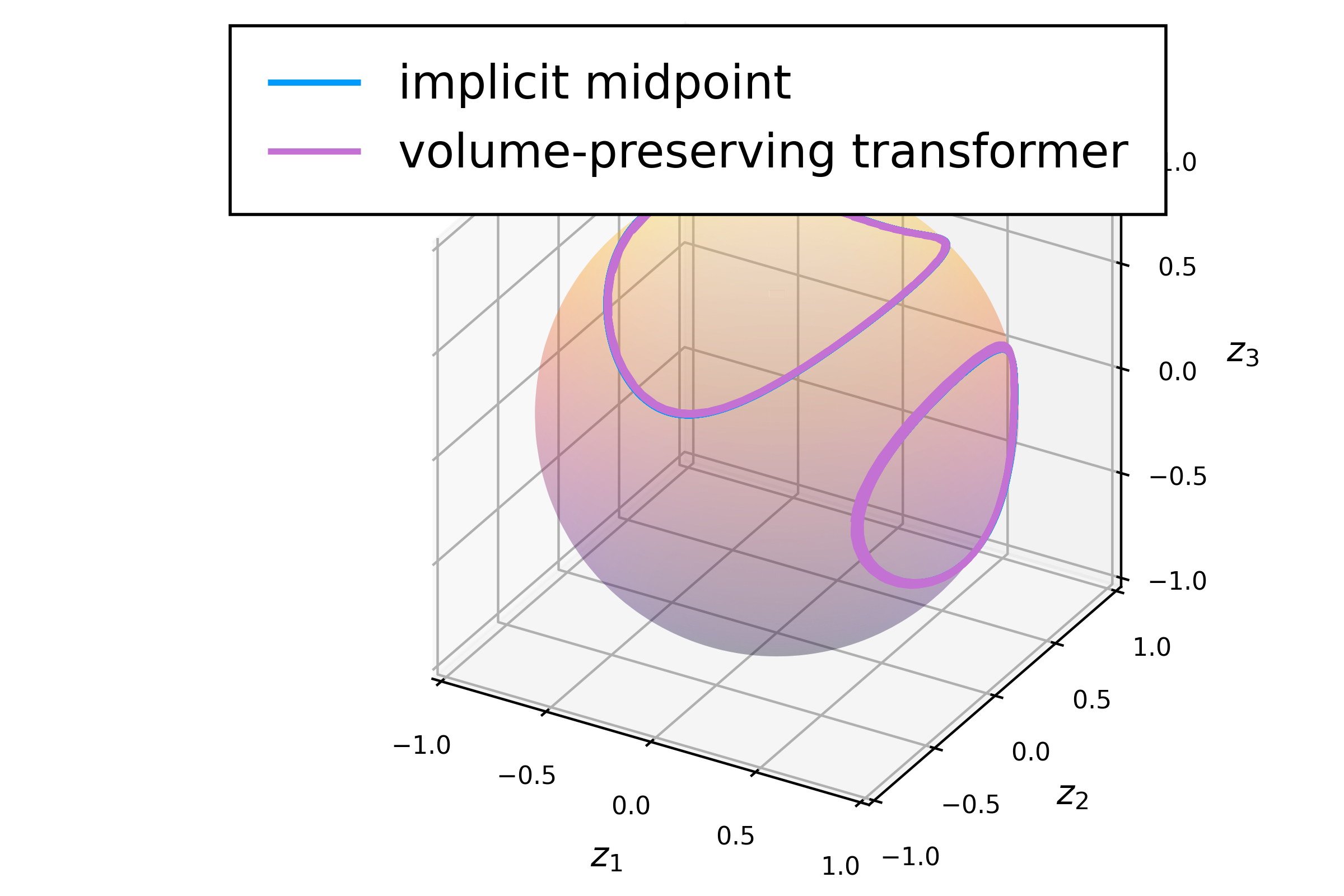}%
\includegraphics[width = .33\textwidth]{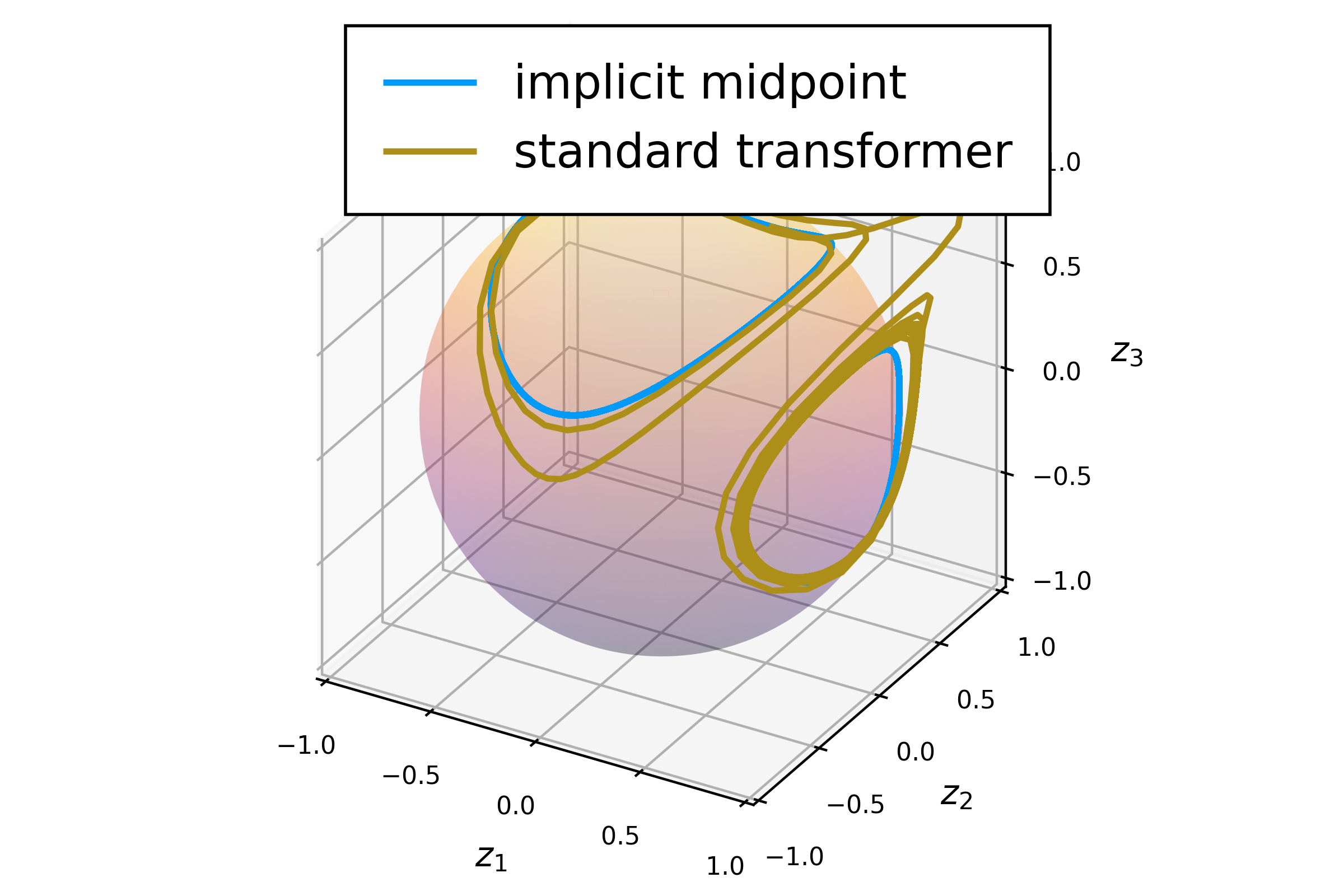}
\caption{Sample trajectories of the rigid body obtained with the three neural networks: volume-preserving feedforward, volume-preserving transformer and the standard transformer, together with the numerical solution (implicit midpoint) for ``trajectory 1" and ``trajectory 4" in \Cref{fig:RigidBodyCurves}. The volume-preserving feedforward neural network is provided with the initial condition (i.e. $z^{(0)}$) and then starts the prediction and the two transformers are provided with the first three time steps ($z^{(1)}$ and $z^{(2)}$ are obtained via implicit midpoint) and then start the prediction. The prediction is made for the time interval $[0, 100]$, i.e. 500 time steps in total.}
\label{fig:Validation3d}
\end{figure}

The corresponding predicted trajectories are shown in \Cref{fig:Validation3d}, where we plot the predicted trajectories for two initial conditions, \(( \sin(1.1), \, 0, \, \cos(1.1) )^T\) and \(( 0, \, \sin(1.1) , \, \cos(1.1) )^T\), for the time interval \([0, 100]\). These initial conditions are also shown in \Cref{fig:RigidBodyCurves} as ``trajectory 1" and ``trajectory 4".

\begin{rmrk} Note that for the two transformers we need to supply three vectors as input as opposed to one for the feedforward neural network. We therefore compute the first two steps with implicit midpoint and then proceed by using the transformer. We could however also use the volume-preserving feedforward neural network to predict those first two steps instead of using implicit midpoint. This is necessary if the differential equation is not known.\end{rmrk}

\begin{figure}
\centering
\includegraphics[width = .7\textwidth]{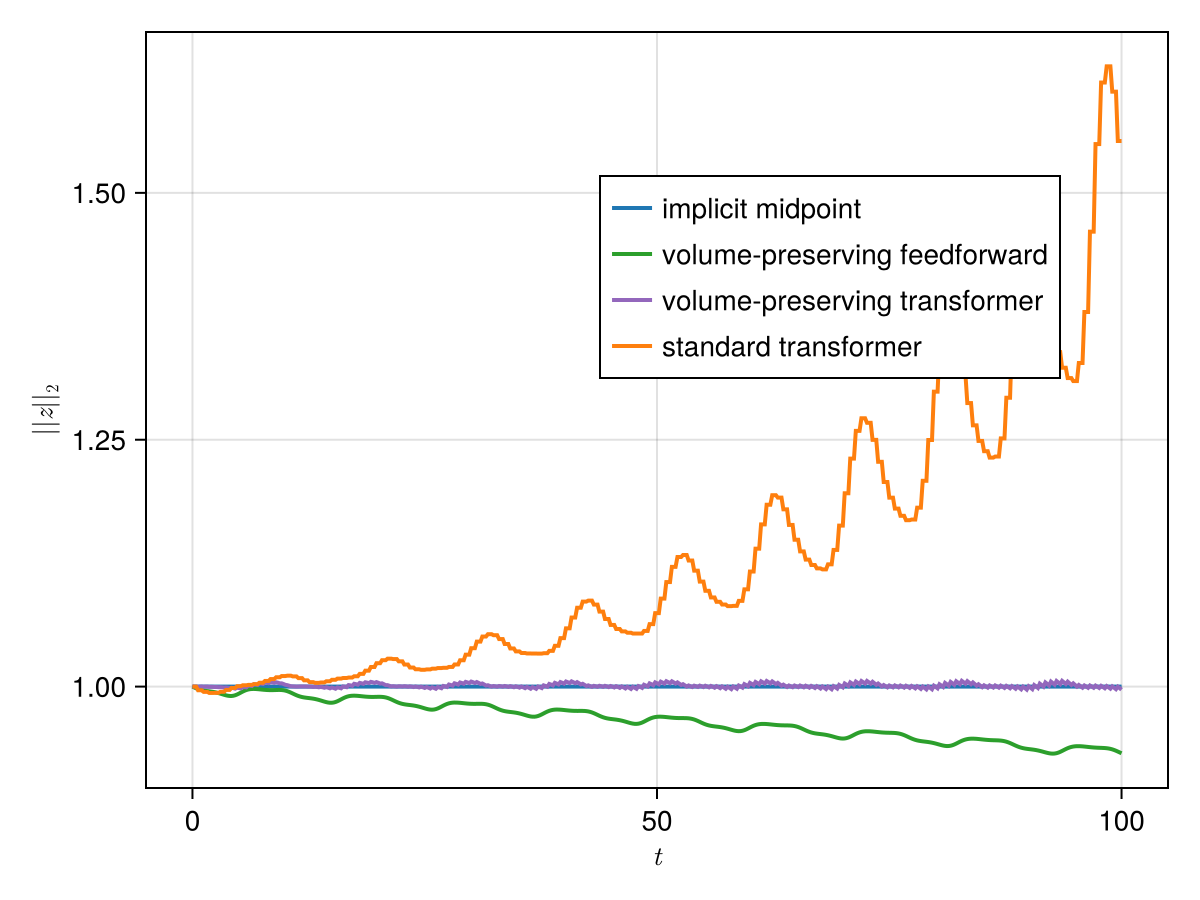}
\caption{Time evolution of invariant $I(z) = \sqrt{z_1^2 + z_2^2 + z_3^2} = ||z||_2$ for ``trajectory 1" up for the time interval $[0, 100]$. We see that for implicit midpoint this invariant is conserved and for the volume-preserving transformer it oscillates around the correct value.}
\label{fig:VPFFvsVPT}
\end{figure}

The standard transformer clearly fails on this task while the volume-preserving feedforward network slowly drifts off. The volume-preserving transformer shows much smaller errors and manages to stay close to the numerical solution.  \Cref{fig:VPFFvsVPT} shows the time evolution of the invariant \(I(z) = ||z||_2\) for implicit midpoint and the three neural network integrators.

In order to get an estimate for the different computational times we perform integration up to time 50000 for all four methods. On CPU we get:

\begin{table}[h]
\centering
\footnotesize\begin{tabulary}{\linewidth}{R L L L L}
\toprule
Method & IM & VPFF & VPT & ST \\
\toprule
Evaluation time & 2.51 seconds & 6.44 seconds & 0.71 seconds & 0.20 seconds \\
\bottomrule
\end{tabulary}

\end{table}

We see that the standard transformer is the fastest, followed by the volume-preserving transformer. The slowest is the volume-preserving feedforward neural network. We attempt to explain these findings:

\begin{itemize}
\item we assume that the standard transformer is faster than the volume-preserving transformer because the softmax can be quicker evaluated than our new activation function \Cref{eq:VolumePreservingActivation},

\item we assume that implicit midpoint is slower than the two transformers because it involves a Newton solver and the neural networks all perform explicit operations,

\item the very poor performance of the volume-preserving feedforward neural network is harder to explain. We suspect that our implementation performs all computations in serial and is therefore slower than the volume-preserving transformer by a factor of three, because we have \(\mathrm{L} = 3\) transformer units. It can furthermore be assumed to be slower by another factor of three because the feedforward neural network only predicts one time step at a time as opposed to three time steps at a time for the two transformer neural networks.

\end{itemize}

Another advantage of all neural network-based integrators over implicit midpoint is that it is easily suitable for parallel computation on GPU because all operations are explicit\footnotemark. The biggest motivation for using neural networks to learn dynamics comes however from non-intrusive reduced order modeling as discussed in the introduction; a traditional integrator like implicit midpoint is simply not suitable for this task as we need to recover dynamics from data.

\footnotetext{The implementation of these architectures in \texttt{GeometricMachineLearning.jl} supports parallel computation on GPU.

}

\subsection{Why Does Regular Attention Fail?}

\label{11790294061382476188}{}

\Cref{fig:Validation3d} shows, that the standard transformer fails to predict the time evolution of the system correctly. The reason behind this could be that it is not sufficiently restrictive, i.e., the matrix which is made up of the three columns in the output of the transformer (see \Cref{eq:StandardTransformerOutput}) does not have full rank (i.e. is not invertible); a property that the volume-preserving transformer has by construction. We observe that the ``trajectory 1" and ``trajectory 4" seem to merge at some point, as if there were some kind of attractor in the system. This is not a property of the physical system and seems to be mitigated if we use volume-preserving architectures.

\subsection{A Note on Parameter-Dependent Equations}

\label{7149635217923694843}{}

In the example presented here, training data was generated by varying the initial condition of the system, specifically

\begin{equation}
\begin{split}\{\varphi^t(z^0_\alpha): {t\in(t_0, t_f], z^0_\alpha\in\mathtt{ics}} \},\end{split}\end{equation}

where \(\varphi^t\) is the flow of the differential equation \(\dot{z} = f(z)\), in particular the rigid body from \Cref{eq:RigidBodyEquations}, \(t_0\) is the initial time, \(t_f\) the final time, and \texttt{ics} denotes the set of initial conditions.

In applications such as \emph{reduced order modeling} [\hyperlinkref{17010868457598550642}{12}–\hyperlinkref{4173415763474909947}{14}], one is often concerned with \emph{parametric differential equations} of the form:

\begin{equation}
\begin{split}\dot{z} = f(z; \mu) \text{ for $\mu\in\mathbb{P}$},\end{split}\end{equation}

where \(\mathbb{P}\) is a set of parameters on which the differential equation depends. In the example of the rigid body, these parameters could be the moments of inertia \((I_1, I_2, I_3)\) and thus equivalent to the parameters \((\mathfrak{a},\mathfrak{b},\mathfrak{c})\) in \Cref{eq:RigidBodyEquations}. A normal feedforward neural network is unable to learn such a parameter-dependent system as it \emph{only sees one point at a time}:

\begin{equation}
\begin{split}\mathcal{NN}_\mathrm{ff}: \mathbb{R}^d\to\mathbb{R}^d.\end{split}\end{equation}

Thus a feedforward neural network can only approximate the flow of a differential equation with fixed parameters as the prediction becomes ambiguous in the case of data coming from solutions for different parameters. A transformer neural network\footnotemark{} on the other hand, is able to describe solutions with different parameters of the system, as it is able to \emph{consider the history of the trajectory up to that point}.

\footnotetext{It should be noted that recurrent neural networks such as LSTMs [\hyperlinkref{3150593183825968864}{5}] are also able to do this. 

}

We also note that one could build parameter-dependent feedforward neural networks as

\begin{equation}
\begin{split}\mathcal{NN}_{ff}:\mathbb{R}^d\times\mathbb{P} \to \mathbb{R}^d.\end{split}\end{equation}

A simple single-layer feedforward neural network

\begin{equation}
\begin{split}    \overline{\mathcal{NN}}_{ff}: x \mapsto \sigma(Ax + b),\end{split}\end{equation}

could be made parameter-dependent by modifying it to

\begin{equation}
\begin{split}    \mathcal{NN}_{ff}: (x, \mu) \mapsto \sigma\Big(A\begin{bmatrix}x \\ \mu \end{bmatrix} + b\Big),\end{split}\end{equation}

for example. This makes it however harder to imbue the neural network with structure-preserving properties and we do not pursue this approach further in this work.

\section{Conclusion and future work}

\label{7826571555278476054}{}

We have introduced a new neural network architecture, referred to as \emph{volume-preserving transformer}, for learning the dynamics of systems described by divergence-free vector fields. This new architecture is based on the classical transformer, but modifies the attention mechanism, such that the resulting integrator is volume-preserving. We have shown that the new network leads to more accurate results than both the classical transformer and volume-preserving feedforward neural networks, when applied to volume-preserving dynamical systems, specifically the rigid body.

Future work will focus on the application of the new architecture to different systems and in particular parameter-dependent equations. A more thorough study of why and when classical attention fails (as was demonstrated in this work) is also desirable. Another interesting objective for future research is the proof of a universal approximation theorem for the volume-preserving feedforward neural network and perhaps even the transformer presented in this work. As was already observed by other authors [\hyperlinkref{2750591267727325161}{10}, \hyperlinkref{2928453333953959728}{21}], efforts should also be directed towards making the learning process more efficient by reducing the number of epochs for small networks by using e.g. a Newton optimizer instead of Adam. Lastly the poor performance of the volume-preserving feedforward neural network, when used has an integrator, has to be investigated further.

\begin{acknowledgement}

\label{15488518252468985777}{}

We would like to thank the Centre International de Rencontres Mathématiques (CIRM) for their excellent hospitality that made this project possible. We further thank all participants of the Cemracs 2023 hackathon with who we had stimulating discussions as well as Léopold Trémant for suggesting the right output format of the transformer to us.

\section*{References}

\label{14148394603223710272}{}

{\raggedright

\hangindent=0.33in {\makebox[{\ifdim0.33in<\dimexpr\width+1ex\relax\dimexpr\width+1ex\relax\else0.33in\fi}][l]{[1]}}\hypertarget{9930491254245683687}{}N.~Baker, F.~Alexander, T.~Bremer, A.~Hagberg, Y.~Kevrekidis, H.~Najm, M.~Parashar, A.~Patra, J.~Sethian, S.~Wild and others. \emph{Workshop report on basic research needs for scientific machine learning: Core technologies for artificial intelligence} (USDOE Office of Science (SC), Washington, DC, 2019).

\hangindent=0.33in {\makebox[{\ifdim0.33in<\dimexpr\width+1ex\relax\dimexpr\width+1ex\relax\else0.33in\fi}][l]{[2]}}\hypertarget{3402864322513845913}{}B.~Chen, K.~Huang, S.~Raghupathi, I.~Chandratreya, Q.~Du and H.~Lipson. \emph{Discovering state variables hidden in experimental data}, arXiv~preprint~arXiv:2112.10755 (2021).

\hangindent=0.33in {\makebox[{\ifdim0.33in<\dimexpr\width+1ex\relax\dimexpr\width+1ex\relax\else0.33in\fi}][l]{[3]}}\hypertarget{8146600536394319917}{}P.~Goyal and P.~Benner. \emph{LQResNet: a deep neural network architecture for learning dynamic processes}, arXiv~preprint~arXiv:2103.02249 (2021).

\hangindent=0.33in {\makebox[{\ifdim0.33in<\dimexpr\width+1ex\relax\dimexpr\width+1ex\relax\else0.33in\fi}][l]{[4]}}\hypertarget{2538783340789034789}{}A.~Vaswani, N.~Shazeer, N.~Parmar, J.~Uszkoreit, L.~Jones, A.~N.~Gomez, L.~Kaiser and I.~Polosukhin. \emph{Attention is all you need}. Advances~in~neural~information~processing~systems \textbf{30} (2017).

\hangindent=0.33in {\makebox[{\ifdim0.33in<\dimexpr\width+1ex\relax\dimexpr\width+1ex\relax\else0.33in\fi}][l]{[5]}}\hypertarget{3150593183825968864}{}S.~Hochreiter and J.~Schmidhuber. \emph{Long short-term memory}. Neural~computation \textbf{9}, 1735–1780 (1997).

\hangindent=0.33in {\makebox[{\ifdim0.33in<\dimexpr\width+1ex\relax\dimexpr\width+1ex\relax\else0.33in\fi}][l]{[6]}}\hypertarget{13058814467303700196}{}A.~Dosovitskiy, L.~Beyer, A.~Kolesnikov, D.~Weissenborn, X.~Zhai, T.~Unterthiner, M.~Dehghani, M.~Minderer, G.~Heigold, S.~Gelly and others. \emph{An image is worth 16x16 words: Transformers for image recognition at scale}, arXiv~preprint~arXiv:2010.11929 (2020).

\hangindent=0.33in {\makebox[{\ifdim0.33in<\dimexpr\width+1ex\relax\dimexpr\width+1ex\relax\else0.33in\fi}][l]{[7]}}\hypertarget{18037804960457623862}{}E.~Hairer, C.~Lubich and G.~Wanner. \emph{Geometric Numerical integration: structure-preserving algorithms for ordinary differential equations} (Springer, Heidelberg, 2006).

\hangindent=0.33in {\makebox[{\ifdim0.33in<\dimexpr\width+1ex\relax\dimexpr\width+1ex\relax\else0.33in\fi}][l]{[8]}}\hypertarget{3006316998437427794}{}V.~I.~Arnold. \emph{Mathematical Methods of Classical Mechanics}. \emph{Graduate Texts in Mathematics} (Springer, New York City, 1978).

\hangindent=0.33in {\makebox[{\ifdim0.33in<\dimexpr\width+1ex\relax\dimexpr\width+1ex\relax\else0.33in\fi}][l]{[9]}}\hypertarget{3603161085352894112}{}M.~Raissi, P.~Perdikaris and G.~E.~Karniadakis. \emph{Physics-informed neural networks: A deep learning framework for solving forward and inverse problems involving nonlinear partial differential equations}. Journal~of~Computational~physics \textbf{378}, 686–707 (2019).

\hangindent=0.33in {\makebox[{\ifdim0.33in<\dimexpr\width+1ex\relax\dimexpr\width+1ex\relax\else0.33in\fi}][l]{[10]}}\hypertarget{2750591267727325161}{}P.~Jin, Z.~Zhang, A.~Zhu, Y.~Tang and G.~E.~Karniadakis. \emph{SympNets: Intrinsic structure-preserving symplectic networks for identifying Hamiltonian systems}. Neural~Networks \textbf{132}, 166–179 (2020).

\hangindent=0.33in {\makebox[{\ifdim0.33in<\dimexpr\width+1ex\relax\dimexpr\width+1ex\relax\else0.33in\fi}][l]{[11]}}\hypertarget{12908279347149791397}{}P.~Buchfink, S.~Glas and B.~Haasdonk. \emph{Symplectic Model Reduction of Hamiltonian Systems on Nonlinear Manifolds and Approximation with Weakly Symplectic Autoencoder}. SIAM~Journal~on~Scientific~Computing \textbf{45}, A289-A311 (2023).

\hangindent=0.33in {\makebox[{\ifdim0.33in<\dimexpr\width+1ex\relax\dimexpr\width+1ex\relax\else0.33in\fi}][l]{[12]}}\hypertarget{17010868457598550642}{}T.~Lassila, A.~Manzoni, A.~Quarteroni and G.~Rozza. \emph{Model order reduction in fluid dynamics: challenges and perspectives}. Reduced~Order~Methods~for~modeling~and~computational~reduction, 235–273 (2014).

\hangindent=0.33in {\makebox[{\ifdim0.33in<\dimexpr\width+1ex\relax\dimexpr\width+1ex\relax\else0.33in\fi}][l]{[13]}}\hypertarget{14434722691705629527}{}K.~Lee and K.~T.~Carlberg. \emph{Model reduction of dynamical systems on nonlinear manifolds using deep convolutional autoencoders}. Journal~of~Computational~Physics \textbf{404}, 108973 (2020).

\hangindent=0.33in {\makebox[{\ifdim0.33in<\dimexpr\width+1ex\relax\dimexpr\width+1ex\relax\else0.33in\fi}][l]{[14]}}\hypertarget{4173415763474909947}{}S.~Fresca, L.~Dede and A.~Manzoni. \emph{A comprehensive deep learning-based approach to reduced order modeling of nonlinear time-dependent parametrized PDEs}. Journal~of~Scientific~Computing \textbf{87}, 1–36 (2021).

\hangindent=0.33in {\makebox[{\ifdim0.33in<\dimexpr\width+1ex\relax\dimexpr\width+1ex\relax\else0.33in\fi}][l]{[15]}}\hypertarget{10955802838160861140}{}S.~Yıldız, P.~Goyal and P.~Benner. \emph{Structure-preserving learning for multi-symplectic PDEs}, arXiv~preprint~arXiv:2409.10432 (2024).

\hangindent=0.33in {\makebox[{\ifdim0.33in<\dimexpr\width+1ex\relax\dimexpr\width+1ex\relax\else0.33in\fi}][l]{[16]}}\hypertarget{4181908863425119879}{}T.~M.~Tyranowski and M.~Kraus. \emph{Symplectic model reduction methods for the Vlasov equation}. Contributions~to~Plasma~Physics (2022).

\hangindent=0.33in {\makebox[{\ifdim0.33in<\dimexpr\width+1ex\relax\dimexpr\width+1ex\relax\else0.33in\fi}][l]{[17]}}\hypertarget{6441307755809283184}{}B.~Brantner and M.~Kraus. \emph{Symplectic autoencoders for model reduction of hamiltonian systems}, arXiv~preprint~arXiv:2312.10004 (2023).

\hangindent=0.33in {\makebox[{\ifdim0.33in<\dimexpr\width+1ex\relax\dimexpr\width+1ex\relax\else0.33in\fi}][l]{[18]}}\hypertarget{11021975310350069440}{}S.~Chaturantabut and D.~C.~Sorensen. \emph{Nonlinear model reduction via discrete empirical interpolation}. SIAM~Journal~on~Scientific~Computing \textbf{32}, 2737–2764 (2010).

\hangindent=0.33in {\makebox[{\ifdim0.33in<\dimexpr\width+1ex\relax\dimexpr\width+1ex\relax\else0.33in\fi}][l]{[19]}}\hypertarget{11086911267780014317}{}A.~Hemmasian and A.~Barati Farimani. \emph{Reduced-order modeling of fluid flows with transformers}. Physics~of~Fluids \textbf{35} (2023).

\hangindent=0.33in {\makebox[{\ifdim0.33in<\dimexpr\width+1ex\relax\dimexpr\width+1ex\relax\else0.33in\fi}][l]{[20]}}\hypertarget{4103674955131354433}{}A.~Solera-Rico, C.~S.~Vila, M.~Gómez, Y.~Wang, A.~Almashjary, S.~Dawson and R.~Vinuesa, \emph{\(\beta\)-Variational autoencoders and transformers for reduced-order modelling of fluid flows}, arXiv~preprint~arXiv:2304.03571 (2023).

\hangindent=0.33in {\makebox[{\ifdim0.33in<\dimexpr\width+1ex\relax\dimexpr\width+1ex\relax\else0.33in\fi}][l]{[21]}}\hypertarget{2928453333953959728}{}J.~Bajārs. \emph{Locally-symplectic neural networks for learning volume-preserving dynamics}. Journal~of~Computational~Physics \textbf{476}, 111911 (2023).

\hangindent=0.33in {\makebox[{\ifdim0.33in<\dimexpr\width+1ex\relax\dimexpr\width+1ex\relax\else0.33in\fi}][l]{[22]}}\hypertarget{5266324132587918940}{}F.~Kang and S.~Zai-Jiu. \emph{Volume-preserving algorithms for source-free dynamical systems}. Numerische~Mathematik \textbf{71}, 451–463 (1995).

\hangindent=0.33in {\makebox[{\ifdim0.33in<\dimexpr\width+1ex\relax\dimexpr\width+1ex\relax\else0.33in\fi}][l]{[23]}}\hypertarget{11823669491487214148}{}R.~Guo, S.~Cao and L.~Chen. \emph{Transformer meets boundary value inverse problems}, arXiv~preprint~arXiv:2209.14977 (2022).

\hangindent=0.33in {\makebox[{\ifdim0.33in<\dimexpr\width+1ex\relax\dimexpr\width+1ex\relax\else0.33in\fi}][l]{[24]}}\hypertarget{17560367509306694129}{}B.~Leimkuhler and S.~Reich. \emph{Simulating hamiltonian dynamics} (Cambridge university press, Cambridge, UK, 2004).

\hangindent=0.33in {\makebox[{\ifdim0.33in<\dimexpr\width+1ex\relax\dimexpr\width+1ex\relax\else0.33in\fi}][l]{[25]}}\hypertarget{7708618804017169771}{}S.~Greydanus, M.~Dzamba and J.~Yosinski. \emph{Hamiltonian neural networks}. Advances~in~neural~information~processing~systems \textbf{32} (2019).

\hangindent=0.33in {\makebox[{\ifdim0.33in<\dimexpr\width+1ex\relax\dimexpr\width+1ex\relax\else0.33in\fi}][l]{[26]}}\hypertarget{12095522061475586953}{}R.~L.~Bishop and S.~I.~Goldberg. \emph{Tensor Analysis on Manifolds} (Dover Publications, Mineola, NY, 1980).

\hangindent=0.33in {\makebox[{\ifdim0.33in<\dimexpr\width+1ex\relax\dimexpr\width+1ex\relax\else0.33in\fi}][l]{[27]}}\hypertarget{4243252928265579806}{}J.~Achiam, S.~Adler, S.~Agarwal, L.~Ahmad, I.~Akkaya, F.~L.~Aleman, D.~Almeida, J.~Altenschmidt, S.~Altman, S.~Anadkat and others. \emph{Gpt-4 technical report}, arXiv~preprint~arXiv:2303.08774 (2023).

\hangindent=0.33in {\makebox[{\ifdim0.33in<\dimexpr\width+1ex\relax\dimexpr\width+1ex\relax\else0.33in\fi}][l]{[28]}}\hypertarget{7305568361813525748}{}A.~Graves and A.~Graves. \emph{Long short-term memory}. Supervised~sequence~labelling~with~recurrent~neural~networks, 37–45 (2012).

\hangindent=0.33in {\makebox[{\ifdim0.33in<\dimexpr\width+1ex\relax\dimexpr\width+1ex\relax\else0.33in\fi}][l]{[29]}}\hypertarget{11194799988828072492}{}D.~E.~Rumelhart, G.~E.~Hinton, R.~J.~Williams and others. \emph{Learning internal representations by error propagation} (1985).

\hangindent=0.33in {\makebox[{\ifdim0.33in<\dimexpr\width+1ex\relax\dimexpr\width+1ex\relax\else0.33in\fi}][l]{[30]}}\hypertarget{5557653725674690011}{}K.~He, X.~Zhang, S.~Ren and J.~Sun. \emph{Deep residual learning for image recognition}. In: \emph{Proceedings of the IEEE conference on computer vision and pattern recognition} (2016); pp.~770–778.

\hangindent=0.33in {\makebox[{\ifdim0.33in<\dimexpr\width+1ex\relax\dimexpr\width+1ex\relax\else0.33in\fi}][l]{[31]}}\hypertarget{9125255746634034294}{}R.~T.~Chen, Y.~Rubanova, J.~Bettencourt and D.~K.~Duvenaud. \emph{Neural ordinary differential equations}. Advances~in~neural~information~processing~systems \textbf{31} (2018).

\hangindent=0.33in {\makebox[{\ifdim0.33in<\dimexpr\width+1ex\relax\dimexpr\width+1ex\relax\else0.33in\fi}][l]{[32]}}\hypertarget{10009329608770522327}{}K.~Feng. \emph{The step-transition operators for multi-step methods of ODE{\textquotesingle}s}. Journal~of~Computational~Mathematics, 193–202 (1998).

\hangindent=0.33in {\makebox[{\ifdim0.33in<\dimexpr\width+1ex\relax\dimexpr\width+1ex\relax\else0.33in\fi}][l]{[33]}}\hypertarget{9298822944544965621}{}E.~Hairer. \emph{Conjugate-symplecticity of linear multistep methods}. Journal~of~Computational~Mathematics, 657–659 (2008).

\hangindent=0.33in {\makebox[{\ifdim0.33in<\dimexpr\width+1ex\relax\dimexpr\width+1ex\relax\else0.33in\fi}][l]{[34]}}\hypertarget{9074255080583970709}{}T.~J.~Bridges and S.~Reich. \emph{Multi-symplectic integrators: numerical schemes for Hamiltonian PDEs that conserve symplecticity}. Physics~Letters~A \textbf{284}, 184–193 (2001).

\hangindent=0.33in {\makebox[{\ifdim0.33in<\dimexpr\width+1ex\relax\dimexpr\width+1ex\relax\else0.33in\fi}][l]{[35]}}\hypertarget{15842744304998993410}{}F.~E.~Cellier and E.~Kofman. \emph{Continuous system simulation} (Springer Science \& Business Media, Heidelberg, 2006).

\hangindent=0.33in {\makebox[{\ifdim0.33in<\dimexpr\width+1ex\relax\dimexpr\width+1ex\relax\else0.33in\fi}][l]{[36]}}\hypertarget{16425831850931499293}{}S.~Lang. \emph{Algebra}. Vol.~211 of \emph{Graduate Texts in Mathematics} (Springer, Heidelberg, 2002).

\hangindent=0.33in {\makebox[{\ifdim0.33in<\dimexpr\width+1ex\relax\dimexpr\width+1ex\relax\else0.33in\fi}][l]{[37]}}\hypertarget{7730444816824140145}{}S.~Gowda, Y.~Ma, A.~Cheli, M.~Gwóźzdź, V.~B.~Shah, A.~Edelman and C.~Rackauckas. \href{https://doi.org/10.1145/3511528.3511535}{\emph{High-Performance Symbolic-Numerics via Multiple Dispatch}}. \href{https://doi.org/10.1145/3511528.3511535}{ACM~Commun.~Comput.~Algebra \textbf{55}, 92–96} (2022).

\hangindent=0.33in {\makebox[{\ifdim0.33in<\dimexpr\width+1ex\relax\dimexpr\width+1ex\relax\else0.33in\fi}][l]{[38]}}\hypertarget{10297842816248433659}{}A.~Schwarzenberg-Czerny. \emph{On Matrix Factorization and Efficient Least Squares Solution}. Astronomy~and~Astrophysics~Supplement \textbf{110} (1995).

\hangindent=0.33in {\makebox[{\ifdim0.33in<\dimexpr\width+1ex\relax\dimexpr\width+1ex\relax\else0.33in\fi}][l]{[39]}}\hypertarget{14370501854957973096}{}M.~Kraus. \href{https://doi.org/10.5281/zenodo.3648325}{\emph{GeometricIntegrators.jl: Geometric Numerical Integration in Julia}}, \href{https://github.com/JuliaGNI/GeometricIntegrators.jl}{\texttt{https://github.com/JuliaGNI/GeometricIntegrators.jl}} (2020).

\hangindent=0.33in {\makebox[{\ifdim0.33in<\dimexpr\width+1ex\relax\dimexpr\width+1ex\relax\else0.33in\fi}][l]{[40]}}\hypertarget{287835226255993214}{}J.~Bezanson, A.~Edelman, S.~Karpinski and V.~B.~Shah. \href{https://doi.org/10.1137/141000671}{\emph{Julia: A fresh approach to numerical computing}}. SIAM~review \textbf{59}, 65–98 (2017).

\hangindent=0.33in {\makebox[{\ifdim0.33in<\dimexpr\width+1ex\relax\dimexpr\width+1ex\relax\else0.33in\fi}][l]{[41]}}\hypertarget{17860356218777098896}{}B.~Brantner and M.~Kraus. \emph{GeometricMachineLearning.jl}, \href{https://github.com/JuliaGNI/GeometricMachineLearning.jl}{\texttt{https://github.com/JuliaGNI/GeometricMachineLearning.jl}} (2020).

\hangindent=0.33in {\makebox[{\ifdim0.33in<\dimexpr\width+1ex\relax\dimexpr\width+1ex\relax\else0.33in\fi}][l]{[42]}}\hypertarget{13147958005838726490}{}NVIDIA. \emph{GeForce RTX 4090}, \href{https://www.nvidia.com/de-de/geforce/graphics-cards/40-series/rtx-4090/}{\texttt{https://www.nvidia.com/de-de/geforce/graphics-cards/40-series/rtx-4090/}} (2022).

\hangindent=0.33in {\makebox[{\ifdim0.33in<\dimexpr\width+1ex\relax\dimexpr\width+1ex\relax\else0.33in\fi}][l]{[43]}}\hypertarget{14686158125302097530}{}T.~Besard, C.~Foket and B.~De Sutter. \emph{Effective Extensible Programming: Unleashing Julia on GPUs}. \href{https://doi.org/10.1109/TPDS.2018.2872064}{IEEE~Transactions~on~Parallel~and~Distributed~Systems} (2018), \href{https://arxiv.org/abs/1712.03112}{arXiv:1712.03112 [cs.PL]}.

\hangindent=0.33in {\makebox[{\ifdim0.33in<\dimexpr\width+1ex\relax\dimexpr\width+1ex\relax\else0.33in\fi}][l]{[44]}}\hypertarget{6288960670092762841}{}D.~P.~Kingma and J.~Ba. \emph{Adam: A method for stochastic optimization}, arXiv~preprint~arXiv:1412.6980 (2014).

\hangindent=0.33in {\makebox[{\ifdim0.33in<\dimexpr\width+1ex\relax\dimexpr\width+1ex\relax\else0.33in\fi}][l]{[45]}}\hypertarget{12925252812113732274}{}I.~Goodfellow, Y.~Bengio and A.~Courville. \emph{Deep learning} (MIT press, Cambridge, MA, 2016).

\hangindent=0.33in {\makebox[{\ifdim0.33in<\dimexpr\width+1ex\relax\dimexpr\width+1ex\relax\else0.33in\fi}][l]{[46]}}\hypertarget{16796601878836832165}{}J.~Nocedal and S.~J.~Wright. \emph{Numerical optimization} (Springer, Heidelberg, 1999).

\hangindent=0.33in {\makebox[{\ifdim0.33in<\dimexpr\width+1ex\relax\dimexpr\width+1ex\relax\else0.33in\fi}][l]{[47]}}\hypertarget{5993914653893793734}{}D.~D.~Holm, T.~Schmah and C.~Stoica. \emph{Geometric mechanics and symmetry: from finite to infinite dimensions}. Vol.~12 (Oxford University Press, Oxford, UK, 2009).

}

\end{acknowledgement}\newpage\appendix

\section{The Rigid Body}

\label{6923796705636635695}{}

\begin{figure}[h]
\centering
\includegraphics[width=.5\textwidth]{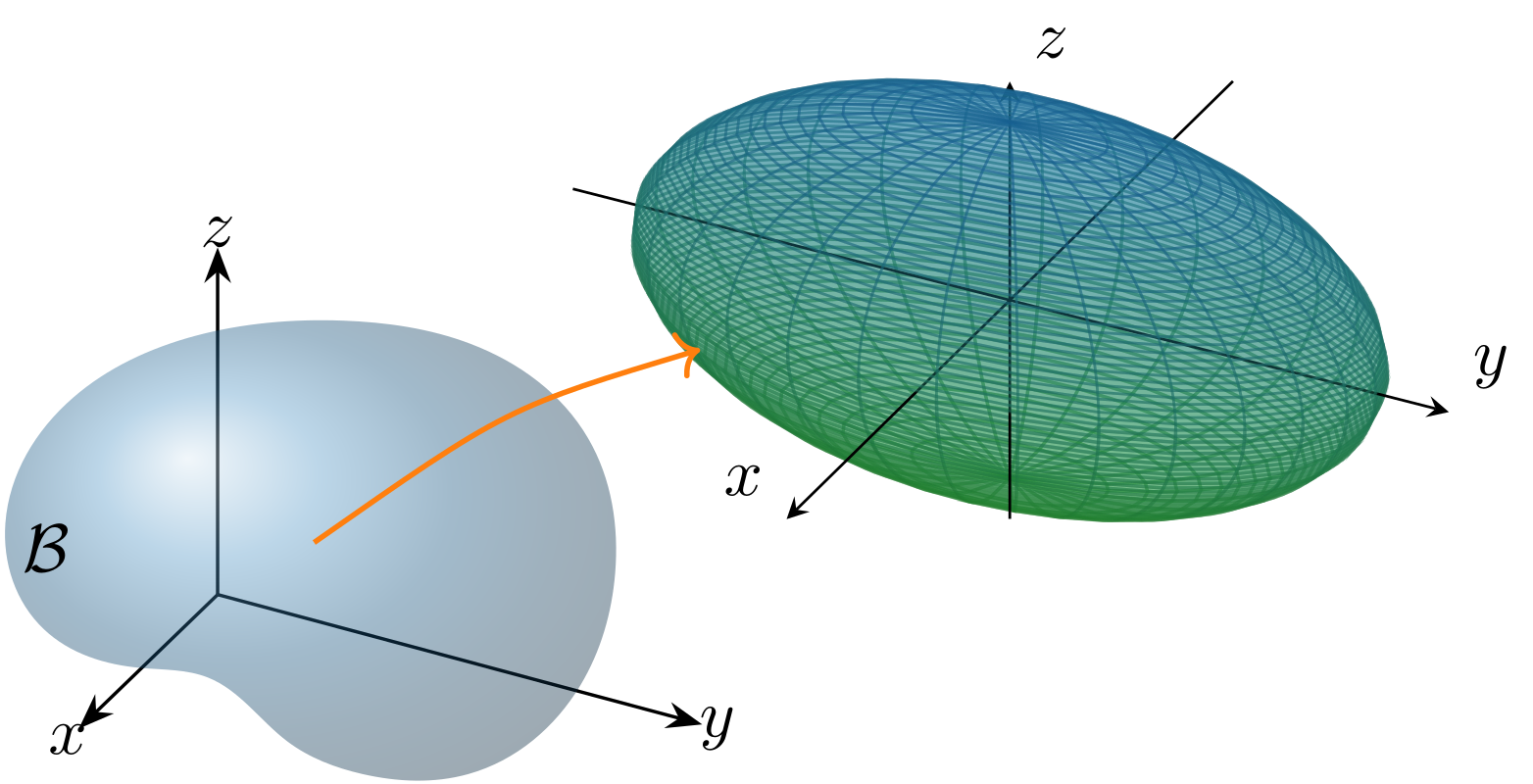}
\caption{Any rigid body fixed at a point (left) can be described through an ellipsoid (right) through $I_1$, $I_2$ and $I_3$.}
\label{fig:RigidBody}
\end{figure}

In this appendix, we sketch the derivation of the rigid body equations used in \Cref{13575308147144999289}. The differential equation for the rigid body [\hyperlinkref{18037804960457623862}{7}, \hyperlinkref{3006316998437427794}{8}] describes the dynamics of a solid object fixed at a point. \Cref{fig:RigidBody} shows an example for a rigid body.  Its dynamics can always be described in terms of an ellipsoid. To see this, let us consider the derivation of the rigid body equations. The motion of a point \((x_1, x_2, x_3)^T\) in the rigid body \(\mathcal{B}\) can be described as follows:

\begin{equation}
\begin{split}v := \frac{d}{dt} \begin{bmatrix} x_1 \\ x_2 \\ x_3 \end{bmatrix}
= \omega \times \begin{bmatrix} x_1 \\ x_2 \\ x_3 \end{bmatrix}
= \begin{pmatrix} \omega_2x_3 - \omega_3x_2 \\ \omega_3x_1 - \omega_1x_3 \\ \omega_1x_2 - \omega_2x_1 \end{pmatrix}
= \begin{pmatrix} 0 & - \omega_3 & \omega_2 \\ \omega_3 & 0 & -\omega_1 \\ -\omega_2 & \omega_1 & 0 \end{pmatrix}\begin{bmatrix} x_1 \\ x_2 \\ x_3 \end{bmatrix},\end{split}\end{equation}

where \(\omega\) is the angular velocity.  The total kinetic energy of \(\mathcal{B}\) is obtained by integrating over the entire volume of the body:

\begin{equation}
\begin{split}T  = \frac{1}{2} \int_\mathcal{B} ||v||^2 dm = \frac{1}{2} \int_\mathcal{B} || \omega \times x ||^2 dm = \frac{1}{2} \omega^T \Theta \omega,        \end{split}\end{equation}

where

\begin{equation}
\begin{split}\Theta_{ij} = \begin{cases}
    \int_\mathcal{B}(x_k^2 + x_\ell^2)dm & \text{ if } i = j, \\  
    -  \int_\mathcal{B}x_ix_jdm & \text{else} ,
\end{cases}\end{split}\end{equation}

where \(dm\) indicates an integral over the mass density of the rigid body, and we further have \(i\neq{}k\neq\ell\neq{}i\) for the first case.

The mathematical description of the motion of a rigid body hence does not require knowledge of the precise shape of the body, but only the coefficients \(\Theta_{ij}\). As \(\Theta\) is a symmetric matrix, we can write the kinetic energy as:

\begin{equation}
\begin{split}T = \frac{1}{2}\omega^T\Theta\omega = \frac{1}{2}\omega^TU^TIU\omega,
\end{split}\label{eq:RigidBodyKineticEnergy}\end{equation}

where \(U\) is an orthonormal matrix that diagonalizes \(\Theta.\) In \Cref{eq:RigidBodyKineticEnergy} we called the eigenvalues of \(\Theta\) \(I = \mathrm{diag} (I_1, I_2, I_3)\).

This shows that it is sufficient to know the eigenvalues of the matrix \(\Theta\) which are called \emph{moments of inertia} and denoted by \(I_k\) for \(k = 1, 2, 3\) to describe the motion of the rigid body (modulo a rotation). From this point of view every rigid body is equivalent to an ellipsoid as indicated in \Cref{fig:RigidBody}.

\subsection{Formulation of the Equations of Motion in the Euler-Poincaré Framework}

\label{670372339511842849}{}

The dynamics of the rigid body can be described through a rotational matrix \(Q(t)\), i.e. each point of the rigid body \(x(0)\in\mathcal{B}\) can be described through \(x(t) = Q(t)x(0)\) where \(Q(t)^TQ(t) = \mathbb{I}\). We can therefore describe the evolution of the rigid body through a differential equation on the Lie group \(G := \{Q\in\mathbb{R}^{d\times{}d}:Q^TQ = \mathbb{I}\}\). The associated tangent vector \(\dot{Q}\in{}T_QG\) can be mapped to the Lie algebra\footnotemark{} \(\mathfrak{g}=T_\mathbb{I}G\) by:

\begin{equation}
\begin{split}W := \dot{Q}Q^T = \begin{pmatrix} 0 & -\omega_3 & \omega_2 \\ \omega_3 & 0 & -\omega_1 \\ -\omega_2 & \omega_1 & 0 \end{pmatrix} \in \mathfrak{g}. 
\end{split}\label{eq:LieAlgebraRepresentation}\end{equation}

\footnotetext{This is the Lie algebra of skew-symmetric matrices: \(\mathfrak{g} = \{W\in\mathbb{R}^{d\times{}d}:W^T = -W\}\).

}

As was indicated in equation \Cref{eq:LieAlgebraRepresentation}, the components of the skew-symmetric matrix \(W\) are equivalent to those of the angular velocity \(\omega\) as can easily be checked:

\begin{equation}
\begin{split}\dot{x} (t)
= \frac{d}{dt}Q(t)x(0) = \dot{Q}(t)x(0) = \dot{Q}(t)Q^T(t)x(t) = W(t)x(t)
= \omega (t) \times x (t) .\end{split}\end{equation}

With this description, the kinetic energy can be written as:

\begin{equation}
\begin{split}T = \frac{1}{2}\mathrm{tr}(WDW^T),
\end{split}\label{eq:KineticEnergyForLieGroup}\end{equation}

where \(D = \mathrm{diag}(d_1, d_2, d_3)\) is a diagonal matrix\footnotemark{} that satisfies \(I_1 = d_2 + d_3\), \(I_2 = d_3 + d_1\) and \(I_3 = d_1 + d_2\).

\footnotetext{This matrix is equivalent to the diagonal entries of the \emph{coefficient of inertia matrix} in [\hyperlinkref{5993914653893793734}{47}].

}

We now write \(z := I^{-1}\omega\) and introduce the following notation\footnotemark:

\footnotetext{Note that the \(\hat{}\) operation used here is different from the hat operation used in \Cref{18309819361304734052}.

}

\begin{equation}
\begin{split}\hat{z} = \widehat{\begin{bmatrix} z_1 \\ z_2 \\ z_3 \end{bmatrix}} = \widehat{\begin{bmatrix} \frac{\omega_1}{I_1} \\ \frac{\omega_2}{I_2} \\ \frac{\omega_3}{I_3} \end{bmatrix}} = \begin{pmatrix} 0 & -\frac{\omega_3}{I_3} & \frac{\omega_2}{I_2} \\ \frac{\omega_3}{I_3} & 0 & -\frac{\omega_1}{I_1} \\ -\frac{\omega_2}{I_2} & \frac{\omega_1}{I_1} & 0 \end{pmatrix}.\end{split}\end{equation}

and obtain via the Euler-Poincaré equations\footnotemark{} for \Cref{eq:KineticEnergyForLieGroup}:

\begin{equation}
\begin{split}\widehat{I\dot{\omega}} = [\widehat{I\omega}, W],\end{split}\end{equation}

\footnotetext{For the Euler-Poincaré equations we have to compute variations of \Cref{eq:KineticEnergyForLieGroup} with respect to \(\delta{}W = \delta(\dot{Q}Q^{-1}) = \dot{\Sigma} - [W,\Sigma]\) where \(\Sigma := \delta{}QQ^{-1}\). For more details on this see [\hyperlinkref{5993914653893793734}{47}].

}

or equivalently:

\begin{equation}
\begin{split}\frac{d}{dt}\begin{bmatrix} z_1 \\  z_2 \\ z_3  \end{bmatrix} 
= \begin{bmatrix} \mathfrak{a}z_2z_3 \\ \mathfrak{b}z_1z_3 \\ \mathfrak{c}z_1z_2 \end{bmatrix}.
\end{split}\label{eq:FinalRigidBodyEquations}\end{equation}

In the above equation, we defined \(\mathfrak{a} := I_3^{-1} - I_2^{-1}\), \(\mathfrak{b} := I_1^{-1} - I_3^{-1}\) and \(\mathfrak{c} := I_2^{-1} - I_1^{-1}\). In all of the examples, we set \(I_1 = 1\), \(I_2 = 2\) and \(I_3 = 2/3\), thus yielding \(\mathfrak{a} = 1\), \(\mathfrak{b} = -1/2\) and \(\mathfrak{c} = -1/2\).

\end{document}